\documentclass[
11pt, reqno]{amsart}

\usepackage{amsthm, amsfonts, amssymb}
\theoremstyle{definition}
\newtheorem{ntn}{Notation}[section]
\newtheorem{dfn}[ntn]{Definition}
\theoremstyle{plain}
\newtheorem{lem}[ntn]{Lemma}
\newtheorem{prp}[ntn]{Proposition}
\newtheorem{thm}[ntn]{Theorem}
\newtheorem{cor}[ntn]{Corollary}

\theoremstyle{remark}
\newtheorem{rem}[ntn]{Remark}
\newtheorem{exa}[ntn]{Example}

\def\floor[#1]{\lfloor #1 \rfloor }

\newcommand{\z}{\mathbb{Z}}

\newcommand{\q}{\mathbb{Q}}
\newcommand{\R}{\mathbb{R}}
\newcommand{\HH}{\mathbb{H}}

\newcommand{\lan}{\langle}
\newcommand{\ran}{\rangle}

\newcommand{\sr}{{\rm sr}(R)}

\newcommand{\GL}{\mathit{{\rm GL}}}

\newcommand{\PGL}{\mathit{{\rm PGL}}}
\newcommand{\SL}{\mathit{{\rm SL}}}
\newcommand{\PSL}{\mathit{{\rm PSL}}}
\newcommand{\jac}{\mathit{{\rm Jac}}}

\newcommand{\SK}{\mathit{{\rm SK}}}

\newcommand{\E}{\mathcal{E}}
\newcommand{\D}{\mathfrak{d}}

\newcommand{\inc}{{\rm inc}}
\newcommand{\id}{{\rm id}}
\newcommand{\tor}{{{\rm Tor}_1^{A}}}
\newcommand{\tors}{{{\rm Tor}_1^{\z}}}
\newcommand{\zzz}{\z[\frac{1}{2}]}

\newcommand{\zzzz}{\z[\frac{1}{n}]}

\newcommand{\s}{\Sigma}
\newcommand{\si}{\sigma}

\newcommand{\arr}{\rightarrow}
\newcommand{\larr}{\longrightarrow}
\newcommand{\harr}{\hookrightarrow}

\newcommand{\se}{\subseteq}

\newcommand{\mt}{\mapsto}
\newcommand{\two}{\twoheadrightarrow}

\newcommand{\fff}{{F^\ast}}

\newcommand{\eee}{{\tilde{E}}}

\newcommand{\rr}{{R^\ast}}

\newcommand{\stabe}{{\rm Stab}}
\newcommand{\diag}{{\rm diag}}
\renewcommand{\ker}{{\rm ker}}
\newcommand{\coker}{{\rm coker}}
\newcommand{\im}{{\rm im}}
\newcommand{\ind}{{\rm ind}}
\newcommand{\class}{{\rm cl}}

\newcommand {\mtx}[4]
{\left(
\begin{array}{cc}
#1 & #2   \\
#3 & #4
\end{array}
\right)}

\newtheoremstyle{athm}
  {}
  {}
  {\itshape}
  {}
  {\scshape}
  {}
  {.5em}
  {\thmnote{#3}}
\theoremstyle{athm}
\newtheorem*{athm}{}

\begin{document}

\title{Third homology of general linear groups}
\author{B. Mirzaii}

\begin{abstract}
The third homology group of $\GL_n(R)$ is studied,
where $R$ is a `ring with many units' with center $Z(R)$.
The main theorem states that if
$K_1(Z(R))\otimes \q \simeq K_1(R)\otimes \q $,
(e.g. $R$ a commutative ring  or  a central simple algebra),
then $H_3(\GL_2(R), \q) \arr H_3(\GL_3(R), \q)$ is injective.
If $R$ is commutative, $\q$ can be replaced by a field $k$ such that
 $1/2 \in k$.
For  an infinite field $R$ (resp. an infinite field $R$ such that $\rr=\rr^2$),
we get a better result that
$H_3(\GL_2(R), \zzz) \arr H_3(\GL_3(R), \zzz)$
(resp. $H_3(\GL_2(R), \z) \arr H_3(\GL_3(R), \z)$)
is injective. As an application we study the third homology
group of $\SL_2(R)$ and the indecomposable part of $K_3(R)$.
\end{abstract}

\maketitle

\section{Introduction}

The Hurewicz theorem relates homotopy groups to homology groups, which
are much easier to calculate. This in turn provides a homomorphism from
the Quillen $K_n$-group of a ring $R$ to the $n$-th integral homology of
stable linear group $\GL(R)$, $h_n: K_n(R) \arr H_n(\GL(R),\z)$. One can also
define Milnor $K$-groups, $K_n^M(R)$, and
when $R$ is commutative there is a canonical map
$K_n^M(R) \arr K_n(R)$ \cite{guin1989}.

One of the approaches to investigate $K$-groups is by means
of the homology stability. Suslin's stability theorem states
that for an infinite field  $F$, the natural map
\[
H_i(\GL_n(F),\z) \arr H_i(\GL(F),\z)
\]
is bijective if $n \ge i$ \cite{sus14}.
Using this result Suslin  constructed a map from
$H_n(\GL_n(F), \z)$ to
$K_n^M(F)$ such that the sequence
\begin{gather*}
H_n(\GL_{n-1}(F), \z) \overset{H_n(\inc)}{\larr} H_n(\GL_n(F), \z)
\larr K_n^M(F) \larr 0
\end{gather*}
is exact. Combining these two results he constructed a map from
$K_n(F)$ to  $K_n^M(F)$ such that the composite homomorphism
\begin{gather*}
K_n^M(F) \arr K_n(F) \arr K_n^M(F)
\end{gather*}
coincides with the multiplication by
$(-1)^{n-1}(n-1)!$ \cite[Sec. 4]{sus14}.

These results have been  generalized by Nesterenko-Suslin \cite{nes-sus}
to  commutative local rings with infinite residue fields,
and by Sah \cite{sah1986} and Guin \cite{guin1989}
to a wider class of rings which we call `rings with many units'.

Except for $n=1, 2$, there is no precise information about
the kernel of $H_n(\inc)$.
In this direction Suslin posed a problem,
which is now referred to as `a conjecture by Suslin'
(see \cite[7.7]{bor-yang}, \cite[4.13]{sah}).
\begin{athm}[{\bf Injectivity Conjecture.}]
For any infinite field $F$ the natural homomorphism
\[
H_{n}(\GL_{n-1}(F), \q) \arr H_n(\GL_n(F), \q)
\]
is injective.
\end{athm}
This conjecture is easy if $n=1,2$. For $n=3$ the conjecture was proved
positively by Sah \cite{sah} and Elbaz-Vincent \cite{elb}.
The case $n=4$ is proved by the author in \cite{mirzaii2005}.
The conjecture is proved in full for number fields by Borel and Yang
\cite{bor-yang}.

When $n=3$, in \cite{elb}
Elbaz-Vincent proves the
conjecture for wider class of commutative rings
(called H1-ring in \cite{elb}). In fact he proves that
for any commutative `ring with many units'
$
H_{3}(\GL_{2}(R), \q) \arr H_3(\GL_3(R), \q)$
is injective. We will generalize  this further,
including some class of non-commutative rings.

The above conjecture says that the kernel of $H_n(\inc)$
is in fact torsion.
Our main goal, in this paper, is to study the map $H_3(\inc)$
in such way that we lose less information on its kernel.
Here is our main result.

\begin{athm}[{\bf Theorem \ref{mir-elb4}.}]
Let $R$ be a ring with many units with center $Z(R)$.
Let $k$ be a field such that $1/2 \in k$.
\par {\rm (i)}
If $K_1(Z(R)) \otimes \q \simeq K_1(R)\otimes \q$,
then $H_3(\GL_2(R), \q) \arr H_3(\GL_3(R), \q)$ is injective.
If $R$ is commutative, then $\q$ can be replaced by $k$.
\par {\rm (ii)} If $R$ is an infinite field or a quaternion algebra
over an infinite field, then $H_3(\GL_2(R),\zzz) \arr H_3(\GL_3(R), \zzz)$
is injective.
\par {\rm (iii)}
Let $R=\R$ or let $R$ be an infinite field such that $\rr=\rr^2$.
Then
$H_3(\GL_2(R), \z) \arr H_3(\GL_3(R), \z)$ is injective.
\par {\rm (iv)}
$H_3(\GL_2(\HH), \z) \arr H_3(\GL_3(\HH), \z)$ is bijective,
where $\HH$ is the ring of quaternion.
\end{athm}


Examples of non-commutative rings with many units which
satisfy the condition $K_1(Z(R)) \otimes \q \simeq K_1(R)\otimes \q$
of (i) in the above theorem are Azumaya algebras over commutative
local rings with infinite residue fields.

As an application we generalize and give an easier
proof of the main theorem of  Sah in \cite[Thm. 3.0]{sah}.
Our proof of the next theorem
avoids the case by  case analysis done in \cite{sah}.

\begin{athm}[{\bf Theorem \ref{sl-in-sl}.}]
Let $R$ be a commutative ring with many units.
Let $k$ be a field such that $1/2 \in k$.
\par {\rm (i)}
The map $H_0(\rr,H_3(\SL_2(R),k)) \arr H_3(\SL(R),k)$ is injective.
\par {\rm (ii)} For an infinite field $R$,
$H_0(\rr,H_3(\SL_2(R),\zzz)) \arr H_3(\SL(R),\zzz)$ is injective.
\par {\rm (iii)} If $R=\R$ or if $R$ is an infinite field such that $\rr=\rr^2$,
then  $H_3(\SL_2(R), \z) \arr H_3(\SL(R), \z)$ is injective.
\par {\rm (iv)}
The map $H_3(\SL_2(\HH), \z) \arr H_3(\SL_3(\HH), \z)$ is bijective.
\end{athm}

We use these results to study the third $K$-group of a field.
Let $K_3(R)^{\ind}=\coker(K_3^M(R) \arr K_3(R))$ be the
indecomposable part of $K_3(R)$. In this article we prove that
if $R$ is an infinite field,
\[
K_3(R)^{\ind} \otimes \zzz \simeq H_0(\rr, H_3(\SL_2(R), \zzz)).
\]
Furthermore if $\rr=\rr^2$
or $R=\R$, then
\[
K_3(R)^{\ind} \simeq H_3(\SL_2(R), \z).
\]

To prove these claims, our general strategy will be the same as in
\cite{sah} and \cite{elb}. We will
introduce some spectral sequences similar to ones in \cite{elb},
smaller but still big enough to do some computation.
The main theorem will come out of analysis of these spectral sequences.

Here we establish some notations.
In this paper by $H_i(G)$ we mean the $i$-th integral homology of the
group $G$. We use the bar resolution to define the homology of a
group \cite[Chap. I, Section 5]{bro}. Define
${\rm \bf{c}}({g}_1, {g}_2,\dots, {g}_n)=
\sum_{\si \in \s_n}
{{\rm sign}(\si)}
[{g}_{\si(1)}| {g}_{\si(2)}|\dots | {g}_{\si(n)}] \in H_n(G)$,
where ${g}_i \in G$ pairwise commute and
$\s_n$ is the symmetric group of degree $n$.
By $\GL_n$ we mean the general linear group $\GL_n(R)$, where
$R$ is a {\it ring with many units}.
By $Z(R)$ we will mean the center of $R$.

Note that $\GL_0$ is the trivial group and
$\GL_1=\rr$. By $\rr^m$ we mean
$\rr \times \cdots\times \rr$ ($m$-times) or, when
$R$ is commutative and $m \ge 2$,
the subgroup of $\rr$, $\{a^m| a \in \rr\}$, depending on the context.
This will not cause any confusion.
The $i$-th factor of $\rr^m=\rr \times \cdots\times \rr$, ($m$-times),
is denoted by $R_i^\ast$.

\subsection*{Acknowledgements}
I would like to thank W. van der Kallen for his interest in this
work and for his valuable comments. The last draft of this article
was written during my stay at mathematics department of Queen's
University Belfast which I was supported by EPSRC R1724PMR. I would
like to thank them for their support and hospitality.

\section{Rings with many units}\label{rings}

The study of `rings with many units' is originated by  W. van der
Kallen in \cite{vdkallen1977}\footnote{This notion is introduced by
W. van der Kallen.}, where he shows that $K_2$ of such
commutative rings behave very much like $K_2$ of fields.
According to \cite{vdkallen1977} in order to have a nice description
of $K_2(R)$ in terms of generators and relations or in order to have
a nice stability property for $K_2(R)$, the ring should have
`enough invertible elements', and the `more invertible elements' the
ring has, the better description of $K_2(R)$ one gets. In this
direction, see \ref{sk1-k2} for a homological proof of a theorem of
Van der Kallen \cite{vdkallen1977}, due to Nesterenko-Suslin
\cite[Cor. 4.3]{nes-sus}.

In \cite{nes-sus} another definition of rings with many units is given, where
they prove very nice homology stability
results for the homology of general linear groups over these rings.
They further prove that when the ring is a local ring with infinite
residue field, the homology stability bound can be very sharp.

In \cite{guin1989} Guin shows that if a ring satisfies both
the definition of Van der Kallen and of Suslin, then
most of the main results of Suslin in \cite{sus14} still are true.
Following \cite{vdkallen1977} and \cite{nes-sus}
we call such rings, {\it rings with many units}.

\begin{dfn}
We say that $R$ is a {\it ring with
many units} if it has the following properties:
\begin{quote}
{\sf (H1) Hypothesis 1.} For any finite number of
surjective linear forms $f_i: R^n \arr R$,
there exist  $v \in  R^n$
such that $f_i(v) \in \rr$, \\
{\sf (H2) Hypothesis 2.} For any $ n \ge 1$, there exist
$n$ elements of the center of $R$ such that the sum
of each nonempty subfamily belongs to $\rr$.
\end{quote}
\end{dfn}

\begin{rem}
(i) { (H1)} implies that the stable range of $R$ is one, $\sr=1$.
\cite[Prop. 1.4]{guin1989}.
\par (ii) {(H1)} implies {(H2)} if $R$ is
 commutative \cite[Prop. 1.3]{guin1989}.
\par (iii)  Property {(H1)} is considered by Van der Kallen
\cite[Sec. 1]{vdkallen1977} and property {(H2)} is studied by
Nesterenko and Suslin \cite[\S 1]{nes-sus}.
\end{rem}

\begin{exa}\label{exa-rings-units}
\par {\rm (i)}
Let $R$ satisfy property (H2). Then
a semilocal ring
$R$ is a ring with many units if and only if
$R/\jac(R)$ is a ring with many units, where
$\jac(R)$ denotes the Jacobson radical of $R$.
\par {\rm (ii)} Product of rings with many units is a ring with
many units.
\par {\rm (iii)} Let $D$ be a finite-dimensional $F$-division algebra,
$F$ an infinite field. Then $M_n(D)$, $n \ge 1$, is a ring with many
units.
\par {\rm (iv)}
Let $F$ be an infinite  field. Then any finite-dimensional
$F$-algebra is a semilocal ring \cite[\S 20]{lam2001}.
Therefore, it is a ring with many units.
\par {\rm (v)} Let $R$ be a commutative semilocal
ring with many units. Then any
Azumaya $R$-algebra is a ring with many units
(see \cite[\S 20]{lam2001}).
\end{exa}

Here we give two known results which are used in the construction of spectral
sequences in the coming section. They show the need for properties
{(H1)} and {(H2)}.

\begin{lem}\label{H1}
Let $R$ satisfy the property {\rm (H1)}.
Let $n \ge 2$ and assume $T_i$, $1\le i \le l$, are finitely many
finite subsets of $R^n$ such that each $T_i$ is a basis of a free
summand of $R^n$ with $k$ elements, where $k \le n-1$.
Then there is a vector $v \in R^n$, such that
$T_i\cup\{v\}$, $1\le i \le l$, is a basis of a free summand of $R^n$.
\end{lem}
\begin{proof}
This is  well-known and easy to prove. We leave the proof to the reader.
\end{proof}

The next result is due to Suslin.

\begin{prp}\label{H2}
Let $R$ satisfy the property {\rm (H2)}.
Let $G_i$ be a subgroup of $\GL_{n_i}$, $i=1,2$, and assume that at
least one of them contains the subgroup of diagonal matrices.
Let $M$ be a submodule of $M_{n_1,n_2}(R)$ such that
$G_1M=M=MG_2$. Then the inclusion
\begin{gather*}
\left(
\begin{array}{cc}
G_1 & 0   \\
0 &G_2
\end{array}
\right)
\arr
\left(
\begin{array}{cc}
G_1 & M   \\
0 &G_2
\end{array}
\right)
\end{gather*}
induces isomorphism in homology with coefficients in $\z$.
\end{prp}
\begin{proof}
See \cite[Thm. 1.9]{sus14}.
\end{proof}

The next proposition is rather well-known. We refer the reader to
\cite[3.2]{guin1989} for the definition of
the Milnor $K$-groups $K_n^M(R)$ of a ring $R$.

\begin{prp}\label{sk1-k2}
Let $R$ be a commutative ring with many units. Then
\par {\rm (i)} $\SK_1(R)=0$.
\par {\rm (ii)}(Van der Kallen \cite{vdkallen1977})
\[
K_2(R)\simeq K_2^M(R)= \rr \otimes_\z \rr/\lan a \otimes
(1-a):a, 1-a \in \rr \ran.
\]
\end{prp}
\begin{proof}
(i) By the homology stability theorem \cite[Thm. 1]{guin1989}
\[
K_1(R)=H_1(\GL(R)) \simeq H_1(\GL_1(R))\simeq \rr.
\]
But we also have $K_1(R) \simeq \rr \times \SK_1(R)$. Thus $\SK_1(R)=0$. \\
(ii) (Nesterenko-Suslin) By easy analysis of the
Lyndon-Hochschild-Serre spectral sequence associated to
\[
1 \arr \SL \arr \GL \arr \rr \arr 1,
\]
using part (i) and the homology stability theorem , one sees that
$K_2(R)\simeq H_2(\GL_2)/H_2(\GL_1)$
(see
\cite[Lem. 4.2]{nes-sus}).
By \cite[Thm. 2]{guin1989} we have
$K_2^M(R)\simeq H_2(\GL_2)/H_2(\GL_1)$.
Therefore $K_2^M(R)\simeq K_2(R)$. For the rest see
\cite[Prop. 3.2.3]{guin1989}.
\end{proof}
~

{\sf In this paper we always assume that $R$ is a ring with many units.}

\section{The spectral sequences}\label{spectral}

Let $C_l(R^n)$ and $D_l(R^n)$ be the free abelian groups with
a basis consisting of $(\lan v_0\ran, \dots, \lan v_l\ran)$ and
$(\lan w_0\ran, \dots, \lan w_l\ran)$ respectively,
where every $\min\{l+1, n\}$
of $v_i \in R^n$ and every $\min\{l+1, 2\}$ of $w_i \in R^n$
are basis of a free direct summand of $R^n$. By
$\lan v_i\ran$ and $\lan w_i\ran$ we mean the submodules
of $R^n$ generated by $v_i$ and $w_i$ respectively.
Let $\partial_0: C_0(R^n) \arr
C_{-1}(R^n):=\z$, $\sum_i n_i(\lan v_i\ran) \mt \sum_i n_i$ and
$\partial_l=\sum_{i=0}^l(-1)^id_i: C_l(R^n) \arr C_{l-1}(R^n)$,
$l\ge 1$, where
\[
d_i((\lan v_0\ran, \dots, \lan v_l\ran))= (\lan
v_0 \ran, \dots,\widehat{\lan v_i \ran}, \dots, \lan v_l\ran).
\]
Define the differential $\tilde{\partial}_l=\sum_{i=0}^l(-1)^i
{\tilde{d}}_i: D_l(R^n)\arr D_{l-1}(R^n)$ similar to $\partial_l$.
By Lemma \ref{H1} it is easy to see that the complexes
\begin{gather*}
C_\ast: \ \ \ \ \
0 \leftarrow C_{-1}(R^n) \leftarrow C_{0}(R^n)  \leftarrow \cdots \leftarrow
C_{l-1}(R^n) \leftarrow  \cdots
\\
\\
D_\ast: \ \ \ \ \
0 \leftarrow D_{-1}(R^n) \leftarrow D_{0}(R^n) \leftarrow
\cdots
\leftarrow D_{l-1}(R^n) \leftarrow  \cdots
\end{gather*}
are exact. Consider $C_{i}(R^n)$ and $D_{i}(R^n)$ as left
$\GL_n$-module in a natural way
and convert this action to the right action
by the definition $m.g:=g^{-1}m$.

Take a free left $\GL_n$-resolution
$P_\ast \arr  \z$ of $\z$ with trivial $\GL_n$-action.
{}From the double complexes $C_{\ast} \otimes _{\GL_n} P_\ast$ and
$D_{\ast} \otimes _{\GL_n} P_\ast$, using Prop. \ref{H2},
we obtain two first quadrant
spectral sequences converging to zero with
\begin{gather*}
\begin{array}{l}
E_{p, q}^1(n)= \begin{cases} H_q(\rr^p \times \GL_{n-p}) &
\text{if $0 \le p \le n$}\\
H_q(\GL_n, C_{p-1}(R^n)) & \text{if $p\ge n+1$}
\end{cases},\\
\\
{\eee}_{p, q}^1(n)=
\begin{cases}
H_q(\rr^{p} \times \GL_{n-p})& \text{if $0 \le p \le 2$}\\
H_q(\GL_n, D_{p-1}(R^n))& \text{if $p \ge 3.$} \end{cases}
\end{array}
\end{gather*}
For $1 \le p \le n$ and $q \ge 0$, $d_{p, q}^1(n)
=\sum_{i=1}^p(-1)^{i+1}H_q(\alpha_{i, p})$, where
\begin{gather*}
\begin{array}{l}
\alpha_{i, p}: \rr^p \times \GL_{n-p}\arr \rr^{p-1} \times
\GL_{n-p+1}, \\
(a_1, \dots, a_p, A) \mt
 (a_1, \dots,\widehat{a_i}, \dots, a_p, \left(
\begin{array}{cc}
a_i & 0          \\
0   & A
\end{array}
\right)).
\end{array}
\end{gather*}
In particular for $0 \le p \le n$,
$d_{p, 0}^1(n)= \begin{cases} {\rm id}_\z & \text{if $p$ is odd}\\
0 & \text{if $p$ is even}\end{cases}$. So $E_{p, 0}^2(n)=0$ for $
p \le n-1$. It is also easy to see that $E_{n, 0}^2(n)=E_{n+1, 0}^2(n)=0$.
See the proof of \cite[Thm. 3.5]{mirzaii2003} for more details.

In this note we will use ${\eee}_{p, q}^i(n)$ and
$E_{p, q}^i(n)$ only for $n=3$,
so from now on by ${\eee}_{p, q}^i$ and $E_{p, q}^i$
we mean ${\eee}_{p, q}^i(3)$ and $E_{p, q}^i(3)$ respectively.
We describe $\eee_{p, q}^1$  for $p=3, 4$. Let
\begin{gather*}
w_1=(\lan e_1 \ran, \lan e_2 \ran, \lan e_3 \ran), \  w_2=(\lan
e_1 \ran, \lan e_2 \ran, \lan e_1+ e_2\ran) \in D_2(R^3)
\end{gather*}
and $u_1,\dots, u_5, u_{6, a} \in D_3(R^3)$, $a, a-1 \in \rr$, where
\begin{gather*}
\begin{array}{ll}
u_1= (\lan e_1 \ran, \lan e_2 \ran,
\lan e_3 \ran, \lan e_1+ e_2+e_3 \ran),&
u_2= (\lan e_1 \ran, \lan e_2 \ran,
\lan e_3 \ran, \lan e_1+ e_2 \ran),\\
u_3= (\lan e_1 \ran, \lan e_2 \ran,
\lan e_3 \ran, \lan e_2+e_3 \ran),&
u_4= (\lan e_1 \ran, \lan e_2 \ran,
\lan e_3 \ran, \lan e_1+e_3 \ran),\\
u_5= (\lan e_1 \ran, \lan e_2 \ran, \lan e_1+e_2 \ran, \lan
e_3\ran),&
u_{6, a}= (\lan e_1 \ran, \lan e_2 \ran,
\lan e_1+e_2 \ran, \\
 &
 \hspace{4 cm}
\lan e_1+ae_2 \ran),
\end{array}
\end{gather*}
(see \cite[Lemma 3.3.3]{guin1989}).
By the Shapiro lemma
\begin{gather*}
\begin{array}{l}
\eee_{3, q}^1=H_q(\stabe_{\GL_3}(w_1)) \oplus H_q(\stabe_{\GL_3}(w_2)),\\
\eee_{4, q}^1=\bigoplus_{j=1}^5 H_q(\stabe_{\GL_3}(u_j)) \oplus
[\bigoplus_{a,a-1 \in \rr}H_q(\stabe_{\GL_3}(u_{6, a}))].
\end{array}
\end{gather*}
So by Prop. \ref{H2} we get
\begin{gather*}
\begin{array}{l}
\eee_{3, q}^1=H_q(\rr^3) \oplus H_q(\rr I_2 \times \rr),\\
\eee_{4, q}^1=H_q(\rr I_3) \oplus H_q(\rr I_2 \times \rr) \oplus
H_q(\rr \times \rr I_2) \oplus H_q(T)\\
\hspace{0.9 cm}
\oplus H_q(\rr I_2 \times \rr)
\oplus [\bigoplus_{a,a-1 \in \rr} H_q(\rr I_2 \times \rr)],
\end{array}
\end{gather*}
where $T=\{(a, b, a) \in R^3: a, b \in \rr\}$. Note that
${\tilde{d}}_{p, q}^1=d_{p, q}^1$ for $ p=1, 2$,
${\tilde{d}}_{3, q}^1|_{H_q(\rr^3)}=d_{3, q}^1$ and
${\tilde{d}}_{3, q}^1|_{H_q(\rr I_2 \times \rr)}=H_q(\inc)$,
where $\inc: \rr I_2 \times \rr \arr \rr^3$.

\begin{lem}\label{elb1}
The group $\eee_{p, 0}^2$ is trivial for $0 \le p \le 5$.
\end{lem}
\begin{proof}
Triviality of $\eee_{p, 0}^2$ is easy for $0 \le p \le 2$. To prove the
triviality of $\eee_{3, 0}^2$, note that $\eee_{2, 0}^1=\z$,
$\eee_{3, 0}^1=\z \oplus \z$ and ${\tilde{d}}_{3, 0}^1((n_1,
n_2))=n_1+n_2$, so if $(n_1, n_2) \in \ker({\tilde{d}}_{3, 0}^1)$,
then $n_2=-n_1$. It is easy to see that  this is contained in
$\im({\tilde{d}}_{4, 0}^1)$. We prove the triviality of $\eee_{5, 0}^2$.
Triviality of $\eee_{4, 0}^2$ is similar but much easier.
This proof is just taken from \cite[Sec. 1.3.3]{elb}.\\
{\bf Triviality of $\eee_{5, 0}^2$}. The proof will be in
four steps;\\
{ \bf Step 1}. The sequence
$0 \arr C_\ast(R^3) \otimes_{\GL_3} \z \arr D_\ast(R^3)\otimes_{\GL_3} \z
\arr Q_\ast(R^3) \otimes_{\GL_3} \z \arr 0$ is exact, where
$Q_\ast(R^3):= D_\ast(R^3)/C_\ast(R^3)$.\\
{\bf Step 2}. The group $H_4(Q_\ast(R^3) \otimes_{\GL_3} \z)$ is trivial.\\
{\bf Step 3}. The map induced in homology by
$ C_\ast(R^3) \otimes_{\GL_3} \z \arr D_\ast(R^3) \otimes_{\GL_3} \z$ is
zero in degree $4$.\\
{\bf Step 4}. The group $\eee_{5, 0}^2$ is trivial.\\
{\bf Proof of step 1}. For
$i \ge -1$, $D_i(R^3)\simeq C_i(R^3)\oplus Q_i(R^3)$.
This decomposition is compatible with the action of
$\GL_3$, so we get an exact sequence of $\GL_3$-modules
\[
0 \arr C_i(R^3) \arr D_i(R^3)\arr Q_i(R^3) \arr 0
\]
which splits as a sequence of $\GL_3$-modules. One can easily
deduce the desired exact sequence from this. Note that this
exact sequence does not split as complexes.\\
{\bf Proof of step 2}. The complex $Q_\ast(R^3)$
induces a spectral sequence
\begin{gather*}
{\hat{E}}_{p, q}^1=
\begin{cases}
0& \text{if $0 \le p \le 2$}\\
H_q(\GL_3, Q_{p-1}(R^3))&
\text{if $p \ge 3$ } \end{cases}
\end{gather*}
which converges to zero. To prove the claim it is sufficient  to prove
that ${\hat{E}}_{5, 0}^2=0$ and to prove this it is
sufficient to prove that
${\hat{E}}_{3, 1}^2=0$. One can see that
${\hat{E}}_{3,1}^1=H_1(\rr I_2 \times \rr)$. If $w=( \lan e_1
\ran, \lan e_2 \ran, \lan e_3 \ran, \lan e_1+e_2 \ran) \in
Q_3(R^3)$, then $H_1(\stabe_{\GL_3}(w)) \simeq H_1(\rr I_2 \times
\rr)$ is a summand of ${\hat{E}}_{4, 1}^1$ and ${\hat{d}}_{4,
1}^1: H_1(\stabe_{\GL_3}(w)) \arr {\hat{E}}_{3,1}^1$ is an
isomorphism. So ${\hat{d}}_{4, 1}^1$ is surjective and therefore
${\hat{E}}_{3, 1}^2=0$.\\
{\bf Proof of step 3}. Consider the following commutative diagram
\begin{gather*}
\begin{array}{ccccc}
 C_5(R^3) \otimes_{\GL_3} \z& \arr &  C_4(R^3) \otimes_{\GL_3} \z&
\arr &  C_3(R^3) \otimes_{\GL_3} \z \\
\Big\downarrow & & \Big\downarrow & & \Big\downarrow  \\
D_5(R^3) \otimes_{\GL_3} \z & \arr & D_4(R^3) \otimes_{\GL_3} \z&
\arr &  D_3(R^3) \otimes_{\GL_3} \z .
\end{array}
\end{gather*}
The generators of $C_4(R^3) \otimes_{\GL_3} \z$ are of the form
$x_{a, b} \otimes 1$, where $x_{a, b}=(\lan e_1 \ran, \lan e_2 \ran,
\lan e_1+ ae_2+ be_3 \ran, \lan e_3 \ran, \lan e_1+e_2+e_3
\ran)$, $a, a-1, b, b-1, a-b \in \rr$ (see \cite[Lemma 3.3.3]{guin1989}). Since
$C_3(R^3) \otimes_{\GL_3} \z=\z$,
$(x_{a, b}-x_{c, d})\otimes 1 \in \ker(\partial_4 \otimes 1)$ and
the elements of this form generate $\ker(\partial_4 \otimes 1)$.
Hence  to prove this step it is sufficient to prove that
$(x_{a, b}-x_{c, d})\otimes 1 \in \im(\tilde{\partial}_5 \otimes 1)$.\\
Set $w_{a}'=(\lan e_1 \ran, \lan e_2 \ran, \lan e_1+ ae_2+ e_3
\ran, \lan e_3 \ran, \lan e_1+e_2 \ran, \lan e_1+ ae_2 \ran) \in
D_5(R^3)$, where $a, a-1, \in \rr$. Let $g$,
$g'$, and $g''$ be the matrices
\begin{gather*}
\left(
\begin{array}{ccc}
0 & a^{-1}&0    \\
-1 & 1+a^{-1}&0\\
0 & 0 & 1
\end{array}
\right),
\qquad
\left(
\begin{array}{ccc}
1 & 0 & -1   \\
0 & 1 &-a\\
0 & 0 & 1
\end{array}
\right),
\qquad
\left(
\begin{array}{ccc}
1 & 0 & 0          \\
0 & a^{-1} & 0\\
0 & 0 & 1
\end{array}
\right),
\end{gather*}
respectively, then
\begin{gather*}
g(\tilde{d}_1(w_{a}'))=\tilde{d}_0(w_{a}'),\
g'(\tilde{d}_3(w_{a}'))=\tilde{d}_2(w_{a}'),\
g''(\tilde{d}_4(w_{a}'))=v_{1}'
\end{gather*}
and so $(\tilde{\partial}_5 \otimes 1)(w_{a}'\otimes 1)=
(v_{1}'-v_{a}')\otimes 1$, where
\[
v_{a}'=(\lan e_1 \ran, \lan e_2 \ran,
\lan e_1+ ae_2+ e_3 \ran, \lan e_3 \ran, \lan e_1+e_2 \ran).
\]
Note that the elements of the form $(gw-w)\otimes 1$ are zero in
$D_\ast \otimes_{\GL_3} \z$. If
\begin{gather*}
\begin{array}{l}
u_{a}'=(\lan e_3 \ran,
\lan e_1+ ae_2+ e_3 \ran, \lan e_1 \ran, \lan e_1+e_2 \ran,
\lan e_1+ ae_2 \ran),\\
u_{a}''=(\lan e_1+ ae_2+ e_3 \ran, \lan e_1 \ran,\lan e_2 \ran,
\lan e_1+e_2 \ran, \lan e_1+ ae_2 \ran),
\end{array}
\end{gather*}
where $a, a-1 \in \rr$, then
\begin{gather*}
\begin{array}{l}
gu_{a}'=(\lan e_3 \ran,
\lan e_1+ ae_2+ e_3 \ran, \lan e_2 \ran, \lan e_1+e_2 \ran,
\lan e_1+ ae_2 \ran),\\
g' u_{a}''=(\lan e_3 \ran, \lan e_1 \ran,\lan e_2 \ran,
\lan e_1+e_2 \ran,\lan e_1+ ae_2 \ran).
\end{array}
\end{gather*}
So if $a,a-1, c,c-1 \in \rr$, then
\[
({\tilde{\partial}}_5 \otimes 1)((z_{a}- z_{c})\otimes 1)=
(t_{c}-t_{a})\otimes 1,
\]
where
\begin{gather*}
\begin{array}{l}
z_{a}=(\lan e_3 \ran,
\lan e_1+ ae_2+ e_3 \ran, \lan e_1 \ran,\lan e_2 \ran, \lan e_1+e_2 \ran,
\lan e_1+ ae_2 \ran),\\
t_{a}=(\lan e_3 \ran, \lan e_1+ ae_2+ e_3 \ran,
\lan e_1 \ran, \lan e_2 \ran, \lan e_1+e_2 \ran).
\end{array}
\end{gather*}
If $g_1$, $g_2$, $g_3$ and $g_4$ are the matrices
\begin{gather*}
\left(
\begin{array}{ccc}
-1 & 0 & 1         \\
-1 & 0 &0\\
\frac{b-1}{1-a} & \frac{1-b}{1-a} & 0
\end{array}
\right),
\left(
\begin{array}{ccc}
0 & -1 & 1          \\
0 & -1 & 0\\
\frac{b-a}{1-a} & \frac{a-b}{1-a} & 0
\end{array}
\right),
\left(
\begin{array}{ccc}
1 & 0 & -1          \\
0 & 1 & -1\\
0 & 0 & \frac{1-b}{b}
\end{array}
\right),
\left(
\begin{array}{ccc}
1 & 0 & 0          \\
0 & 1 & 0\\
0 & 0 & \frac{1}{b}
\end{array}
\right),
\end{gather*}
respectively,
then
\begin{gather*}
\begin{array}{ll}
g_1(\tilde{d}_0(y_{a, b}))=t_{\frac{1}{1-b}},   &
g_2(\tilde{d}_1(y_{a, b}))=t_{\frac{-a}{b-a}},  \\
g_3(\tilde{d}_3(y_{a, b}))=v_{\frac{a-b}{1-b}}',&
g_4(\tilde{d}_3(y_{a, b}))=v_{a}',
\end{array}
\end{gather*}
where
\begin{gather*}
y_{a, b}=(\lan e_1 \ran, \lan e_2 \ran,
\lan e_1+ ae_2+ be_3 \ran, \lan e_3 \ran, \lan e_1+e_2+e_3
\ran,  \lan e_1+e_2 \ran).
\end{gather*}
(Here by $\frac{r}{s} \in \rr$ we mean $s^{-1}r$.) By an easy computation
\begin{gather*}
({\tilde{\partial}}_5 \otimes 1)(y_{a, b}\otimes 1)=
t_{\frac{1}{1-b}}\otimes 1 -
t_{\frac{-a}{b-a}} \otimes 1 +
v_{1}' \otimes 1 
\\
\hspace{3 cm}
-v_{\frac{a-b}{1-b}}' \otimes 1 +
v_{a}' \otimes 1 - x_{a, b} \otimes 1.
\end{gather*}
Now it is easy to see that $(x_{a, b}-x_{c, d})\otimes 1 \in
({\tilde{\partial}}_5 \otimes 1)(D_5(R^3)\otimes_{\GL_3} \z)$. This
completes the proof of step 3.\\
{\bf Proof of Step 4}. From the homology long exact
sequence of the short
exact sequence obtained in the first step, we get the exact sequence
\[
H_4(C_\ast(R^3)\otimes_{\GL_3} \z) \arr H_4(D_\ast(R^3)\otimes_{\GL_3} \z)
\arr H_4(Q_\ast(R^3)\otimes_{\GL_3} \z).
\]
By steps 2 and 3,  $H_4(D_\ast(R^3)\otimes_{\GL_3} \z)=0$,
but $\eee_{5, 0}^2=H_4(D_\ast(R^3)\otimes_{\GL_3} \z)$. This completes the
proof of the triviality of $\eee_{5, 0}^2$.
\end{proof}

\begin{lem}\label{elb2}
The group $\eee_{p, 1}^2$ is trivial for $0 \le p \le 4$.
\end{lem}
\begin{proof}
Triviality of $\eee_{p, 1}^2$, $p=0,1$, is a result
of Lemma \ref{elb1} and the fact that the spectral sequence
converges to zero (one can also prove this directly).
If $(a_0, b_0, c_0) \in \ker(\tilde{d}_{2, 1}^1)$, $a_0, b_0, c_0 \in H_1(\rr)$,
then $a_0=b_0$. It is easy to see that this
element is contained in  $\im(\tilde{d}_{3, 1}^1)$.
Let $x=(x_1, \dots, x_5, (x_{6, a})) \in \eee_{4,1}^1$, where
$x_2=(a_2, a_2, b_2) $, $x_3=(a_3, b_3, b_3)$, $x_4=(a_4, b_4, a_4)$,
$x_5=(a_5, a_5, b_5)$, $a_i, b_i \in H_1(\rr)$. By a direct
calculation ${\tilde{d}}_{4, 1}(x)=( p_1, p_2)$, where
\begin{gather*}
\begin{array}{l}
p_1=-(a_2, a_2, b_2)
-(a_3, b_3, b_3)
+(b_4, a_4, a_4)
+(a_5, a_5, b_5),\\
p_2=(a_2, a_2, b_2)
+(b_3, b_3, a_3)
-(a_4, a_4, b_4)
-(a_5, a_5, b_5).
\end{array}
\end{gather*}
If $y=((a_0, b_0, c_0),(d_0, d_0, e_0)) \in \ker(\tilde{d}_{3, 1}^1)$, $a_0, b_0,
c_0, d_0 ,e_0 \in H_1(\rr)$, then
$b_0+d_0=a_0-b_0+c_0+e_0=0$. Let $x_2'=(-b_0, -b_0, -c_0)$,
$x_3'=(-a_0+b_0, 0, 0)$ and set
$x'=(0, x_2', x_3', 0, 0, 0)\in \eee_{4,1}^1$, then
$y={\tilde{d}}_{4, 1}(x')$.

To prove the triviality of
$\eee_{4,1}^2$; let $x \in \ker({\tilde{d}}_{4, 1})$ and set
\begin{gather*}
\begin{array}{ll}
w_1=(\lan e_1 \ran, \lan e_2 \ran, \lan e_1+e_2 \ran,
\lan e_3 \ran,\lan e_1+ ae_2 \ran),\\
w_2=(\lan e_1 \ran, \lan e_2 \ran, \lan e_3\ran,
\lan e_1+e_3 \ran,\lan e_1+ be_3 \ran),\\
w_3=(\lan e_1 \ran, \lan e_2 \ran, \lan e_3 \ran,
\lan e_1+e_2+e_3\ran, \lan e_2+ e_3 \ran),\\
w_{4, a}=(\lan e_1 \ran, \lan e_2 \ran, \lan e_3\ran,
\lan e_1+e_2\ran, \lan e_1+ ae_2 \ran),\\
w_5=(\lan e_1 \ran, \lan e_2 \ran, \lan e_3 \ran,
\lan e_1+e_2+e_3\ran, \lan e_1+ ae_2+be_3 \ran),
\end{array}
\end{gather*}
where $a,a-1, b, b-1, a-b \in \rr$, $b$ fixed.
The groups  $T_i=H_1(\stabe_{\GL_3}(w_i))$, $i=1,2,3,5$ and
$T_4=\bigoplus_{a,a-1 \in \rr} H_1(\stabe_{\GL_3}(w_{4, a}))$
are summands of $\eee_{5, 1}^1$. Note that
$T_1=H_1(\rr I_2 \times \rr)$, $T_2=H_1(T)$,
$T_3=T_5=H_1(\rr I_3)$ and
$T_4=\bigoplus_{a,a-1 \in \rr}H_1(\rr I_2 \times \rr)$.
The restriction of ${\tilde{d}}_{5, 1}^1$
on these summands is as follows,
\begin{gather*}
\begin{array}{l}
{\tilde{d}}_{5, 1}^1|_{T_1}((c_1, c_1, d_1))=(0, (c_1, c_1, d_1), 0, 0,
(c_1, c_1, d_1), -(c_1, c_1, d_1)),\\
{\tilde{d}}_{5, 1}^1|_{T_2}((c_2, d_2, c_2))=(0, 0,(d_2, c_2, c_2)
, (c_2, d_2, c_2),0 , -(c_2, c_2, d_2)),\\
{\tilde{d}}_{5, 1}^1|_{T_3}((c_3, c_3, c_3))
=((c_3, c_3, c_3), (c_3, c_3, c_3),-(c_3, c_3, c_3),0 ,0 ,0),\\
{\tilde{d}}_{5, 1}^1|_{T_{4, a}}((c_4, c_4, d_4))=
(0, 0 , 0 , 0 , 0 ,(c_4, c_4, d_4)),\\
{\tilde{d}}_{5, 1}^1|_{T_5}= \id_{H_1(\rr I_3)}.
\end{array}
\end{gather*}
Let $z_1=(a_5, a_5, b_5) \in T_1$ and $z_2=(a_4, b_4, a_4) \in T_2$. Then
$x-{\tilde{d}}_{5, 1}^1(z_1+z_2)=(x_1', x_2', x_3', 0, 0, (x_{6, a}'))$,
so we can assume that $x_4=x_5=0$.
An easy calculation shows that $a_2=b_2=- a_3=-b_3$.
If $z_3=(a_2, a_2, a_2) \in T_3$,
then $x-{\tilde{d}}_{5, 1}(z_3)=(x_1', 0 , 0,0 ,0 ,(x_{6, a}'))$.
Again we can assume that $x_2=x_3=0$.
If $z_4=(x_{6, a}) \in T_4$, then
$x-{\tilde{d}}_{5, 1}^1(z_4)=(x_1', 0 , 0 , 0 , 0 , 0)$.
Once more we can assume that $x_{6, a}=0$.
These reduce $x$ to an element of the form $(x_1, 0 , 0 , 0 , 0 ,0)$.
If $x_1 \in T_5$, then
${\tilde{d}}_{5, 1}^1(x_1)= (x_1, 0 , 0 , 0 , 0 , 0)$.
This completes the triviality of $\eee_{4, 1}^2$.
\end{proof}

\begin{lem}\label{elb3}
The group
$\eee_{p, 2}^2$ is trivial for $0 \le p \le 3$.
\end{lem}
\begin{proof}
Triviality of $\eee_{0, 2}^2$ and $\eee_{1, 2}^2$
is a result of \ref{elb1}, \ref{elb2} and the
fact that the spectral sequence converges to zero.
Let
\begin{gather*}
\begin{array}{l}
\eee_{1, 2}^1=H_2(\rr \times \GL_2)=
H_2(\rr) \oplus H_2(\GL_2) \oplus H_1(\rr) \otimes
H_1(\GL_2),\\
\eee_{2, 2}^1=H_2(\rr^3)=\bigoplus_{i=1}^6 T_i,\\
\eee_{3, 2}^1=H_2(\rr^3) \oplus H_2(\rr I_2 \times \rr)=
\bigoplus_{i=1}^9 T_i,
\end{array}
\end{gather*}
where
\begin{gather*}
\begin{array}{ll}
T_i=H_2(R_i^\ast)\ {\rm for}\  i=1, 2, 3, &
T_4=H_1(R_1^\ast)\otimes H_1(R_2^\ast),\\
T_5=H_1(R_1^\ast)\otimes H_1(R_3^\ast),   &
T_6=H_1(R_2^\ast)\otimes H_1(R_3^\ast), \\
T_7=H_2(\rr I_2), &
T_8=H_2(I_2 \times \rr),\\
T_9=H_1(\rr I_2)\otimes H_1(I_2 \times \rr). &
\end{array}
\end{gather*}
If $y=(y_1, y_2, y_3, \sum r \otimes s, \sum t \otimes u, \sum v \otimes w)
\in \eee_{2, 2}^1$ and
\begin{gather*}
\hspace{-0.5 cm}
x=(x_1, x_2, x_3, \sum a\otimes b, \sum c\otimes d,
\sum e\otimes f, x_7,x_8,
\sum g\otimes h) \in \eee_{3, 2}^1,
\end{gather*}
$a, b, \dots ,h ,r , \dots, w \in H_1(\rr)$,
then ${\tilde{d}}_{2, 2}^1(y)=(h_1, h_2, h_3)$, where
\begin{gather*}
\begin{array}{l}
h_1=-y_1+y_2,\\
h_3= -\sum s \otimes \diag(1, r) - \sum r \otimes \diag(1, s) \\
\hspace{0.9 cm}
-\sum t \otimes \diag(1, u)+\sum v \otimes \diag(1, w)
\end{array}
\end{gather*}
and ${\tilde{d}}_{3, 2}^1(x)=(z_i)_{1 \le i \le 6}$, where
\begin{gather*}
\begin{array}{l}
z_1=z_2=x_2+x_7, \ \
z_3= x_1+x_3-x_2+x_8,\\
z_4=\sum a\otimes b- \sum c\otimes d+ \sum e\otimes f,\\
z_5=-\sum b \otimes a-\sum a\otimes b+ \sum c\otimes d+ \sum g\otimes h,\\
z_6=-\sum d \otimes c+\sum f\otimes e+ \sum e\otimes f+ \sum g\otimes h.
\end{array}
\end{gather*}
If $y \in \ker({\tilde{d}}_{2, 2}^1)$, then $y_1=y_2$ and $h_3=0$. By
the isomorphism
$H_1(\rr) \otimes H_1(\GL_1) \simeq H_1(\rr) \otimes H_1(\GL_2)$ and
the triviality of $h_3$ we have
\[
-\sum s \otimes r - \sum r \otimes s-\sum t \otimes u+\sum v \otimes w=0.
\]
If
\[
z=(y_1, y_1, y_3, 0, \sum t \otimes u, \sum r \otimes s +
\sum t \otimes u, 0, 0, 0) \in \eee_{3, 2}^1,
\]
then $y=\tilde{d}_{3, 2}^1(z)$ and therefore $\eee_{2, 2}^2=0$.

Let ${\tilde{d}}_{3, 2}^1(x)=0$. Consider the summands
$S_2= H_2(\stabe_{\GL_3}(u_2))=H_2(\rr I_2 \times \rr)$
and $S_3=H_2(\stabe_{\GL_3}(u_3))=H_2(\rr \times \rr I_2)$
of $\eee_{4,2}^1$. Then $S_i \simeq H_2(\rr) \oplus H_2(\rr) \oplus
H_1(\rr) \otimes H_1(\rr)$ and  by a direct calculation
\begin{gather*}
\begin{array}{l}
{\tilde{d}}_{4, 2}^1|_{S_2}((y_1, y_2, s \otimes t))=
(-y_1, -y_1, -y_2, 0, -s \otimes t,-s \otimes t,y_1, y_2,
s \otimes t),\\
{\tilde{d}}_{4, 2}^1|_{S_3}((q_1, q_2, p \otimes q))\!=\!
(-q_1, -q_2,-q_2, -p \otimes q,- p \otimes q,0,q_2,q_1,-q \otimes p).
\end{array}
\end{gather*}
Choose $z_2'=(-x_2, -x_3, -\sum e \otimes f) \in S_2$ and
$z_3'=(x_3 + x_8, 0, -\sum a \otimes b) \in S_3$.
Then $x={\tilde{d}}_{4, 2}^1(z_2'+z_3')$ and therefore
$\eee_{3,2}^2=0$.
\end{proof}

\begin{lem}\label{mir-elb1}
The groups $\eee_{0, 3}^2$, $\eee_{1, 3}^2$ and $\eee_{0, 4}^3$
are trivial.
\end{lem}
\begin{proof}
This follows from \ref{elb1}, \ref{elb2} and \ref{elb3} and
the fact that the spectral sequence converges to zero.
\end{proof}

\begin{cor}\label{mir-elb2}
{\rm (i)} The complex
\[
H_2(\rr^3\times \GL_{0}) \overset{d_{3,2}^1}{\arr}
H_2(\rr^2\times \GL_{1}) \overset{d_{2,2}^1}{\arr}
H_2(\rr \times \GL_{2}) \overset{d_{1,2}^1}{\arr}
H_2(\GL_{3}) \arr 0
\]
is exact, where
$d_{3,2}^1=H_2(\alpha_{1,3})-H_2(\alpha_{2,3})+H_2(\alpha_{3,3})$,
$d_{2,2}^1=H_2(\alpha_{1,2})-H_2(\alpha_{2,2})$ and
$d_{1,2}^1=H_2(\inc)$.
\par {\rm (ii)} The complex
\[
H_3(\rr^2\times \GL_{1}) \overset{d_{2,3}^1}{\larr}
H_3(\rr \times \GL_{2}) \overset{d_{1,3}^1}{\larr}
H_3(\GL_{3}) \arr 0
\]
is exact, where $d_{2,3}^1=H_3(\alpha_{1,2})-H_3(\alpha_{2,2})$ and
$d_{1,3}^1=H_3(\inc)$.
\end{cor}
\begin{proof}
The case (i) follows from the proof of Lemma \ref{elb3} and (ii) follows
from Lemma \ref{mir-elb1}.
\end{proof}

%

\begin{lem} \label{indecomposable}
The groups $E_{0, 4}^3$, $E_{5, 0}^3$ are trivial.
\end{lem}
\begin{proof}
Using \ref{mir-elb2}, one sees that $E_{p, q}^2$-terms are of the form
\begin{gather*}
\begin{array}{cccccccc}
E_{0,4}^2 & \ast &           &           &       &           &       &\\
0         & 0    & E_{2,3}^2 & \ast      & \ast  & \ast      &       & \\
0         & 0    &  0        & \ast      & \ast  & \ast      &       & \\
0         & 0    &  0        & E_{3,1}^2 & \ast  & \ast      & \ast  &\\
0         & 0    &  0        &  0        & 0     & E_{5,0}^2 & \ast. &
\end{array}
\end{gather*}
{}From this description we get $E_{3,1}^3\simeq E_{3,1}^\infty=0$.
So we obtain the exact sequence
\begin{gather*}
0 \arr E_{5,0}^3 \arr E_{5,0}^2
\overset{d_{5,0}^2}{\arr} E_{3,1}^2 \arr 0.
\end{gather*}
The map of spectral sequences $E_{p,q} \arr \eee_{p,q}$ induces
the following commutative diagram
\[
\begin{array}{ccccc}
E_{3,3}^1 & \overset{d_{3,3}^1}{\larr} &   E_{2,3}^1 &
\overset{d_{2,3}^1}{\larr} & E_{1,3}^1 \\
\Big\downarrow\vcenter{%
\rlap{$\scriptstyle{}$}} & &
\Big\downarrow\vcenter{%
\rlap{$\scriptstyle{}$}}
& & \Big\downarrow   \\
\eee_{3,3}^1 & \overset{{\tilde{d}}_{3,3}^1}{\larr} &
\eee_{2,3}^1 & \overset{{\tilde{d}}_{2,3}^1}{\larr} & \eee_{1,3}^1.
\end{array}
\]
Since $E_{p,q}^1=\eee_{p,q}^1$
for $p= 0, 1,2$, the diagram induces the surjective
map $E_{2,3}^2 \two \eee_{2,3}^2$. Now look at the commutative
diagram
\[
\begin{array}{ccc}
 E_{2,3}^2 &  \overset{d_{2,3}^2}{\larr} & E_{0,4}^2 \\
\Big\downarrow\vcenter{%
\rlap{$\scriptstyle{}$}}
& & \Big\downarrow   \\
\eee_{2,3}^2 & \overset{{\tilde{d}}_{2,3}^2}{\larr} & \eee_{0,4}^2.
\end{array}
\]
From the definitions of the spectral sequences
\[
E_{0,4}^2=\eee_{0,4}^2 = H_4(\GL_3)/\im H_4(\rr \times \GL_2).
\]
By Lemma \ref{mir-elb1}, ${\tilde{d}}_{2,3}^2$ is surjective, so the
surjectivity of $d_{2,3}^2$ follows from the commutativity of the
diagram and the surjectivity of the left-hand column map. Therefore
$E_{0,4}^3=0$.

Using this it is easy to see that $E_{5,0}^3 \simeq E_{5,0}^\infty$. Since
the spectral sequence converges to zero, we have $E_{5,0}^3=0$.
\end{proof}

Following \cite[Section 3]{yag} we define;

\begin{dfn}
Let $F$ be an infinite field. We call
\[
\wp^n(F)_\class := H(C_{n+2}(F^n)_{\GL_n} \arr C_{n+1}(F^n)_{\GL_n}
\arr C_n(F^n)_{\GL_n})
\]
the $n$-th classical Bloch group.
\end{dfn}

\begin{prp}\label{divisible}
Let $F$ be an infinite field.
We have an isomorphism $\wp ^3(F)_\class \simeq \fff$.
In particular if $F$ is algebraically closed, then
$\wp^3(F)_\class$ is divisible.
\end{prp}
\begin{proof}
In the proof of \ref{indecomposable} we obtained the exact sequence
\begin{gather*}
0 \arr E_{5,0}^3 \arr E_{5,0}^2
\overset{d_{5,0}^2}{\arr} E_{3,1}^2 \arr 0.
\end{gather*}
By \ref{indecomposable}, $E_{5,0}^3=0$. By the above definition
$E_{5,0}^2=\wp ^3(F)_\class$. It is also easy to
see that $E_{3,1}^2=H_1(\fff)$.
This proves the first part of the proposition.
The second part follows from  the fact that
for an algebraically closed field $F$, $\fff$ is
divisible.
\end{proof}

\begin{rem}\label{conjn=3}
{}From \ref{divisible} and the existence of a surjective map
$\wp^3(F)_\class \arr \wp^3(F)$ \cite[Prop. 3.11]{yag}
we deduce that $\wp^3(F)$ is divisible.
See \cite[2.7]{yag} for the definition of $\wp^3(F)$. This gives a
positive answer to conjecture $0.2$ in \cite{yag} for $n=3$.
\end{rem}

\section{K\"unneth Theorem for $H_3(\fff \times\fff) $}

Let $F$ be an infinite field.
The K\"unneth
theorem for $H_3(\mu_F \times\mu_F)$ finds the following form
\begin{gather*}
0 \arr H_3(\mu_F) \oplus H_3(\mu_F) \arr H_3(\mu_F \times \mu_F)
\arr \tors( \mu_F,\mu_F) \arr 0.
\end{gather*}
Clearly $H_3(\mu_F)\oplus H_3(\mu_F) \arr H_3(\mu_F \times \mu_F)$ 
is the map  $\alpha:=H_3(i_1)+H_3(i_2)$,
where $i_l: \mu_F \arr
\mu_F \times \mu_F$ is the usual injection, $l=1,2$. Let
\[
\beta: H_3(p_1) \oplus H_3(p_2): H_3(\mu_F \times \mu_F)
\arr H_3(\mu_F) \oplus H_{3}(\mu_F),
\]
where $p_l: \mu_F \times\mu_F \arr \mu_F$
is the usual projection, $l=1,2$. Since $\beta\circ\alpha=\id$,
the above exact sequence
splits canonically. Thus we have the canonical decomposition
\begin{gather*}
H_3(\mu_F \times \mu_F)=H_3(\mu_F) \oplus H_3(\mu_F)
\oplus \tors(\mu_F,\mu_F).
\end{gather*}
We construct a splitting map
$\tors(\mu_F,\mu_F)\arr H_3(\mu_F \times \mu_F)$.
The elements of the group  $\tors(\mu_F,\mu_F)
=\tors(H_1(\mu_F),H_1(\mu_F))$
are of the form $\lan \xi, n, \xi \ran=\lan [\xi], n, [\xi] \ran$,
$\xi$ is an element of order $n$ in $\fff$
\cite[Chap. V, Section 6]{mac} (note that $\tors(\mu_F,\mu_F)\simeq \mu_F$). 
It is easy to see that
$\partial_2(\sum_{i=1}^{n}[\xi| \xi^i])= n[\xi]$ in ${(B_1)}_{\mu_F}$.
For the definition of
$\partial_2$ and $B_\ast$ see \cite[Chap. IV, Sec. 5]{mac}.
By \cite[Chap. V, Prop. 10.6]{mac} a map
$\phi: \tors(H_1(\mu_F),H_1(\mu_F))
\arr H_3({(B_\ast)}_{\mu_F}\otimes {(B_\ast)}_{\mu_F})$
can be defined as
\[
a:=\lan [\xi] , n, [\xi] \ran \mt
[\xi] \otimes \sum_{i=1}^{n}[\xi| \xi^i]  +
\sum_{i=1}^{n}[\xi| \xi^i]\otimes [\xi].
\]
Considering the isomorphism
${(B_\ast)}_{\mu_F}\otimes {(B_\ast)}_{\mu_F} \simeq
{(B_\ast)}_{\mu_F \times \mu_F}$
we have $\phi(a)=\chi(\xi) \in H_3(\mu_F \times \mu_F)$, where
\begin{gather*}
\chi(\xi):= \sum_{i=1}^{n}\big([(\xi,1)|(1,\xi)|(1,\xi^i)]-
[(1,\xi)|(\xi,1)|(1,\xi^i)]
+ [(1,\xi)|(1,\xi^i)|(\xi,1)]\\
\hspace{1.7 cm}
 + [(\xi,1)|(\xi^i,1)|(1,\xi)] -
[(\xi,1)|(1,\xi)|(\xi^i,1)]+
[(1,\xi)|(\xi,1)|(\xi^i,1)]\big).
\end{gather*}
Consider the following commutative diagram
\begin{displaymath}
\begin{array}{ccccccccc}
\! 0  & \!\!\!\!\! \arr &
\!\!\!\!\!\!\!\!\! H_3(\mu_F) \oplus H_{3}(\mu_F) &
\!\!\!\! \arr  & \!\!\!
H_3(\mu_F \times \mu_F) & \!\!\!\!\! \arr &
\!\!\!\! \tors(\mu_F,\mu_F)  & \!\!\!\! \arr & \!\!\!\! 0 \\
& &\Big\downarrow & & \Big\downarrow & & \Big\downarrow & &\\
0  & \!\!\!\!\! \arr & \!\!\!
\bigoplus_{i+j=3} H_i(\fff) \otimes H_j(\fff)
& \!\!\!\! \arr &  \!\!\!\!\! H_3(\fff \times \fff)
& \!\!\!\!\!\arr & \!\!\!\!\! \tors(\fff,\fff) &
\!\!\!\!\!\! \arr & \!\!\!\!\!\! 0.
\end{array}
\end{displaymath}
Since $\tors(\mu_F,\mu_F) \simeq \tors(\fff,\fff)$,
we see that the second horizontal exact sequence in the above diagram
splits canonically. So we proved the following proposition.

\begin{prp}\label{canonical}
Let $F$ be an infinite field.
Then we have the canonical decomposition
\[
H_3(\fff \times \fff)= \bigoplus_{i+j=3} H_i(\fff) \otimes H_j(\fff)
\oplus \tors(\fff,\fff),
\]
where a splitting map
$\tors(\fff,\fff)=\tors(\mu_F,\mu_F)
\arr H_3(\fff \times \fff)$
 is defined by $\lan [\xi] , n, [\xi] \ran \mt \chi(\xi)$.
\end{prp}


\section{The injectivity theorem}

\begin{lem}\label{H-1}
Let $K_1(Z(R))\otimes \zzzz \overset{\theta}{\simeq} K_1(R)\otimes \zzzz$ 
induced by the usual inclusion $Z(R) \arr R$.
Then for all $i \ge 1$, $H_i(Z(R),\zzzz) \simeq H_i(K_1(R), \zzzz)$.
\end{lem}
\begin{proof}
Since the map  $\theta$ is an isomorphism in the localized category of
$\zzzz$-modules, it induces an isomorphism on the group homology 
in this category.
\end{proof}

\begin{exa}\label{k1-zr}
(i) If $R$ is commutative, then $K_1(Z(R))= K_1(R)$.
\par (ii)
Let $R$ be a (finite dimensional) division $F$-algebra of rank
$[R:F]=n^2$. Note that $F=Z(R)$. Then
%
$K_1(F) \otimes \z[\frac{1}{n}] \simeq K_1(R)\otimes \z[\frac{1}{n}]$.
This is also true if $R$ is an  Azumaya $S$-algebra, $S$ a commutative
local ring \cite[Cor. 2.3]{hazrat2006}.

These are the examples one should keep in mind in the rest of
this section.
\end{exa}

Let $A$ be a commutative ring with trivial $\GL_3$-action.
Let $P_\ast \arr A$ be a
free left $A[\GL_3]$-resolution of $A$
with trivial $\GL_3$-action. Consider the complex
\begin{gather*}
D_\ast': \ \
0  \leftarrow   D_0'(R^3)  \leftarrow  D_1'(R^3)
\leftarrow  \cdots \leftarrow  D_l'(R^3) \leftarrow  \cdots,
\end{gather*}
where $D_i'(R^3):=D_i(R^3)\otimes A$. The double complex
$D_\ast' \otimes_{\GL_3} P_\ast$ induces a  first quadrant spectral
sequence ${\E}_{p, q}^1 \Rightarrow H_{p+q}(\GL_3, A)$, where
${\E}_{p, q}^1=\eee_{p+1, q}^1(3)\otimes A$ and
${\D}_{p, q}^1=\tilde{d}_{p+1, q}^1\otimes \id_A$.

\begin{lem}\label{mir-elb3}
The groups
$\E_{3, 0}^2$, $\E_{4, 0}^2$, $\E_{2, 1}^2$, $\E_{3, 1}^2$,
$\E_{1, 2}^2$ and $\E_{2, 2}^2$ are trivial.
\end{lem}
\begin{proof}
This follows from the above spectral sequence and Lemmas
\ref{elb1}, \ref{elb2}, \ref{elb3}.
\end{proof}

\begin{thm}\label{mir-elb4}
Let $Z(R)$ be the center of $R$. Let $k$ be a field such that $1/2 \in k$.
\par {\rm (i)}
If $K_1(Z(R)) \otimes \q \simeq K_1(R)\otimes \q$, then
$H_3(\GL_2, \q) \arr H_3(\GL_3, \q)$ is injective.
If $R$ is commutative, then $\q$ can be replaced by $k$.
\par {\rm (ii)} If $R$ is an infinite field or a quaternion algebra
over an infinite field, then $H_3(\GL_2,\zzz) \arr H_3(\GL_3, \zzz)$
is injective.
\par {\rm (iii)}
Let $R=\R$ or let $R$ be an infinite field such that $\rr=\rr^2$. Then
$H_3(\GL_2) \arr H_3(\GL_3)$ is injective.
\par {\rm (iv)}
$H_3(\GL_2(\HH)) \arr H_3(\GL_3(\HH))$ is bijective.
\end{thm}
\begin{proof}
Let $A=\z$, $\zzz$, $\q$ or $k$ (depending on parts (i),...,(iv)).
By Lemma \ref{mir-elb3},
${\E}_{0, 3}^2\simeq {\E}_{0,3}^\infty \simeq H_3(\GL_3, A)$, so to
prove the theorem it is sufficient to prove that $H_3(\GL_{2},A)$ is
a summand of ${\E}_{0, 3}^2$. To prove this it is sufficient to define a
map $\varphi: H_3(\rr \times \GL_2, A) \arr H_3(\GL_2, A)$ such that
$\varphi|_{H_3(\GL_{2}, A)}$ is the identity map and
${\D}_{1,3}^1(H_3(\rr^2 \times \GL_1, A)) \se \ker(\varphi)$.

We have the canonical decomposition
$H_3(\rr \times \GL_{2}, A)= \bigoplus_{i=0}^4 S_i$,
where
\begin{gather*}
\begin{array}{l}
S_i=H_i(\rr, A) \otimes H_{3-i}(\GL_{2},A),\ \ 0 \le i \le 3 \\
S_4=\tor(H_1(\rr, A), H_1(\GL_2, A)).
\end{array}
\end{gather*}
In case of (i) this follows from the K\"unneth theorem and the fact that
$S_4=0$. In other cases it follows again from the K\"unneth theorem and  an
argument  in the line of the previous section.
Note that for parts (ii), (iii) and (iv), the splitting map is
\begin{gather*}
S_4 \simeq \tors(\mu_{Z(R)},\mu_{Z(R)})\otimes A
 \overset{\phi}{\arr} H_3(\rr \times \rr, A)
\overset{q_\ast}{\arr} H_3(\rr \times \GL_{2}, A),
\end{gather*}
where $\phi$ can be  defined as in the previous section, and
\[
q:\rr \times \rr \arr \rr \times \GL_{2}, \ \ (a, b) \mt (a, \diag(b, 1)).
\]
Define $\varphi|_{S_0}: S_0 \arr H_3(\GL_2, A)$ the identity map,
\[
\varphi|_{S_2}: S_2 \simeq H_2(\rr, A) \otimes
H_{1}(\GL_{1}, A) \arr H_3(\rr \times \GL_1, A) \arr H_3(\GL_2, A)
\]
the shuffle product, $\varphi|_{S_3}: S_3 \arr H_3(\GL_2, A)$ the map induced
by $\rr \arr \GL_2$, $a \mt \diag(a, 1)$,
and $\varphi|_{S_4}: S_4 \arr H_3(\GL_2, A)$ the composition
\[
S_4 \overset{\phi}{\arr}H_3(\rr \times \rr, A)
\overset{\inc_\ast}{\arr} H_3(\GL_2, A).
\]
By homology stability theorem \cite[Thm. 1]{guin1989} and
a theorem of Dennis \cite[Cor. 8]{dennis1976} (see also
\cite[Thm. 1]{arlettaz1990})
we have the  decomposition
\begin{gather*}
H_2(\GL_2)=H_2(K_1(R)) \oplus K_2(R).
\end{gather*}
So using \ref{H-1} we have $S_1=S_1' \oplus S_1''$, where
\begin{gather*}
\begin{array}{l}
S_1'=H_1(\rr, A) \otimes H_2(Z(\rr), A),\\
S_1''=H_1(\rr, A) \otimes K_2(R) \otimes A.
\end{array}
\end{gather*}
Define $\varphi|_{S_1'}: S_1' \arr H_3(\GL_2, A)$ the shuffle product and define
the map $\varphi|_{S_1''}: S_1'' \arr H_3(\GL_2, A)$ as the composition
\begin{gather*}
\hspace{-1.8 cm}
H_1(Z(\rr), A) \otimes K_2(R) \otimes A \overset{f}{\arr}
H_1(Z(\rr), A) \otimes H_2(\GL_2, A)\\
\hspace{5.5 cm}
\overset{g}{\arr} H_3(Z(\rr) \times \GL_2, A)
\overset{h}{\arr} H_3(\GL_2, A),
\end{gather*}
where $f=\frac{1}{2}\lambda$, $\lambda$ the natural map
\begin{gather*}
\lambda: K_2(R) \otimes A = H_2(E(R), A) \arr H_2(\GL(R), A) \simeq H_2(\GL_2, A),
\end{gather*}
$g$ is the shuffle product and $h$ is induced by the map
\[
Z(\rr) \times \GL_2 \arr \GL_2, \ \
(a, B)\mt aB.
\]
By Proposition \ref{canonical} we have
$H_3(\rr^2 \times \GL_1, A)=\bigoplus_{i=0}^8 T_i$, where
\begin{gather*}
\begin{array}{ll}
T_0=H_3(\GL_1, A),\\
T_1=\bigoplus_{i=1}^3 H_i(R_1^\ast, A)\otimes H_{3-i}(\GL_1,A), \\
T_2=\bigoplus_{i=1}^3 H_i(R_2^\ast , A)\otimes H_{3-i}(\GL_1, A), \\
T_3=H_1(R_1^\ast, A)\otimes H_1(R_2^\ast , A)\otimes H_1(\GL_1, A),\\
T_4=\tor(H_1(R_1^\ast, A), H_1(R_2^\ast , A)),\\
T_5= \tor(H_1(R_1^\ast, A), H_1(\GL_1 , A)), \\
T_6=\tor(H_1(R_2^\ast, A), H_1(\GL_1 , A)),\\
T_7=H_1(R_1^\ast, A)\otimes H_2(R_2^\ast , A),\\
T_8=H_2(R_1^\ast, A)\otimes H_1(R_2^\ast , A).
\end{array}
\end{gather*}
Note that here $R_i^\ast=\rr$, $i=1, 2$, is the $i$-th summand of $\rr^2=\rr\times \rr$.
We know that ${\D}_{1,3}^1=\si_1 -\si_2$, where $\si_i=H_3(\alpha_{i, 2})$.
It is not difficult to see that
${\D}_{1,3}^1(T_0 \oplus T_1\oplus T_2\oplus T_7 \oplus T_8)
\se \ker(\varphi)$. Here one should use the isomorphism
$H_1(\GL_1, A) \simeq H_1(\GL_2, A)$.
Now $(\si_1-\si_2)(T_4) \se S_4$,
$\si_1(T_5) \se S_0$ and $\si_2(T_5) \se S_4$,
$\si_1(T_6) \se S_4$ and $\si_2(T_6) \se S_0$.
With this description one can see that
${\D}_{1,3}^1(T_4 \oplus T_5\oplus T_6) \se \ker(\varphi)$.
To finish the proof of the claim we have to prove that
${\D}_{1,3}^1(T_3) \se \ker(\varphi)$.
Let $x=a \otimes b \otimes c \in T_3$. By \ref{H-1}
we may assume that
$a, b, c \in Z(\rr)$. Then
\begin{gather*}
{\D}_{1,3}^1(x)= -b \otimes {\rm \bf{c}}(\diag(a,1), \diag(1, c))
-a \otimes {\rm \bf{c}}(\diag(b,1), \diag(1, c)) \in S_1 \\
\hspace{0.7 cm}
=(-b \otimes {\rm \bf{c}} (a, c)- a \otimes {\rm \bf{c}}(b, c),
b \otimes[a, c]+a \otimes [b, c]) \in S_1' \oplus S_1'',
\end{gather*}
where
\begin{gather*}
[a, c]:={\rm {\bf c}}(\diag(a, 1, a^{-1}), \diag(b, b^{-1},1))
\in H_2(E(R),A)\\
={\rm {\bf c}}(\diag(a, 1), \diag(b, b^{-1}))\ \ \ \in H_2(\GL_2,A).
\end{gather*}
Thus
\begin{gather*}
\begin{array}{l}
\varphi({\D}_{1,3}^1(x))=
-{\rm \bf{c}}(\diag(b, 1), \diag(1, a), \diag(1, c))\\
\hspace{2.2 cm}
-{\rm \bf{c}}(\diag(a, 1), \diag(1, b), \diag(1, c))\\
\hspace{2.2 cm}
+\frac{1}{2}{\rm \bf{c}}(\diag(b, b), \diag(a, 1), \diag(c, c^{-1}))\\
\hspace{2.2 cm}
+\frac{1}{2}{\rm \bf{c}}(\diag(a, a), \diag(b, 1), \diag(c, c^{-1})).
\end{array}
\end{gather*}
Set $p:=\diag(p,1), \overline{q}:=\diag(1,q),
p\overline{q}\overline{r}:=
{\rm \bf{c}}(\diag(p, 1), \diag(1, q), \diag(1, r))$, etc. Conjugation by
$\left(
\begin{array}{cc}
0 &  1         \\
1 & 0
\end{array}
\right)$
induces the equality $p\overline{q}\overline{r}=\overline{p}qr$ and
it is easy to see that $pqr=-qpr$ and
$\overline{p^{-1}}qr=-\overline{p}qr$. With these notations and
the above relations we have
\begin{gather*}
\begin{array}{l}
\varphi({\D}_{1,3}^1(x))=-b\overline{ac}-a\overline{bc}+
\frac{1}{2}(b ac+ ba\overline{c^{-1}}+\overline{b}ac+
\overline{b}a\overline{c^{-1}}) \\
\hspace{2.2 cm}
+\frac{1}{2}(abc+ab\overline{c^{-1}}
+\overline{a}bc+\overline{a}b\overline{c^{-1}})=0.
\end{array}
\end{gather*}
This proves that $H_3(\GL_2, A)$ is a summand of ${\E}_{0, 3}^2$.
This proves (i) and (ii).

The proof of (iii) is almost the same as the proof of (i),
only we need to modify the definition of
the map $f$. If $\rr=\rr^2$, $f$ should be induced by the map
\begin{gather*}
K_2(R)= K_2^M(R)\arr H_2(\GL_2), \ \
\{a, b\} \mt
{\rm {\bf c}}(\diag(\sqrt{a}, 1), \diag(b, b^{-1})).
\end{gather*}
Note that if $R$ is commutative and $\rr=\rr^2$, then
$K_2^M(R)$ is uniquely $2$-divisible \cite[1.2]{bass-tate}, so in this case
$f$ is well-defined.

Now Let $R=\R$. It is well-known that
$K_2^M(\R)=\lan \{-1, -1\} \ran \oplus K_2^M(\R)^\circ$, where
$\lan \{-1, -1\} \ran$ is a group of order 2 generated by $\{-1, -1\}$
and $K_2^M(\R)^\circ$ is a uniquely divisible group. In fact every element of
$K_2^M(\R)$ can be uniquely written as $m\{-1, -1\}+\sum \{a_i, b_i\}$,
$a_i, b_i > 0$ and $m=0$ or $1$. Now we define the map
$K_2^M(\R) \arr H_2(\GL_2(\R))$ by
$\{-1, -1\} \mt 0$ and $\{a, b\} \mt
{\rm {\bf c}}(\diag(\sqrt{a}, 1), \diag(b, b^{-1})) $ for $a, b>0$.

For the proof of (iv) we should mention that
$\R^{>0}= K_1^M(\R)^\circ\simeq K_1(\HH)$
and $K_2^M(\R)^\circ\simeq K_2(\HH)$ \cite[p. 188]{sah}.
Since $K_2(\HH)$ and $H_2(\R^{> 0})$ are uniquely divisible,
the proof of injectivity is similar
to the above approach. Surjectivity  follows
from \cite[Thm. 2]{guin1989} and the fact that
$K_n^M(\HH)$ are trivial for $n \ge 2$ \cite[Remark B.15]{sah1986}.
\end{proof}

\begin{cor}\label{mir-sus}
Let $Z(R)$ be the center of $R$. Let $k$ be a field such that $1/2 \in k$.
\par {\rm (i)}
If $K_1(Z(R)) \otimes \q \simeq K_1(R)\otimes \q$,
then we have the exact sequence
 \[
0 \arr H_3(\GL_2, \q) \arr H_3(\GL_3, \q) \arr
K_3^M(R)\otimes \q \arr 0.
\]
If $R$ is commutative, then $\q$ can be replaced by $k$.
\par {\rm (ii)}
If $R$ is an infinite field or a quaternion algebra with an infinite center,
then we have the split exact sequence
 \[
0 \arr H_3(\GL_2, \zzz) \arr H_3(\GL_3, \zzz) \arr
K_3^M(R)\otimes \zzz \arr 0.
\]
\par {\rm (iii)} Let $R$ be an infinite field
 such that $\rr=\rr^2$. Then
we have the split exact sequence
\[
0 \arr H_3(\GL_2) \arr H_3(\GL_3) \arr K_3^M(R) \arr 0.
\]
\par {\rm (iv)}
We have the (non-split) exact sequence
\[
0 \arr H_3(\GL_2(\R)) \arr H_3(\GL_3(\R)) \arr K_3^M(\R) \arr 0.
\]
\end{cor}
\begin{proof}
The exactness in all cases follows from  \ref{mir-elb4} and the
following exact sequence \cite[Thm. 2]{guin1989}
\[
H_3(\GL_2) \arr H_3(\GL_3) \arr K_3^M(R) \arr 0.
\]
If $R$ is commutative, we have a natural map $K_3^M(R)\arr K_3(R)$
such that the composition
\[
K_3^M(R)\arr K_3(R)\arr H_3(\GL_3) \arr K_3^M(R)
\]
coincides with the multiplication by $2$ \cite[Prop. 4.1.1]{guin1989}.
Now splitting maps can be constructed easily.
%
\end{proof}


\begin{rem}
\par (i) Let $R=M_m(D)$, where $D$ is a finite dimensional division
$F$-algebra. Then  $\GL_n(R)\simeq \GL_{mn}(D)$. So by stability
theorem and \cite[Thm. 2]{guin1989}, $K_i^M(R)=0$ for $m \ge 2$
and $i \ge 2$.
\par (ii) It seems that it is not known whether for a finite
dimensional division $F$-algebra $D$,
$H_2(\GL_1(D), \q) \arr H_2(\GL_2(D), \q)$ is injective.
The only case that is known to us is when $D=\HH$.
This follows from  applying the K\"unneth theorem to
$\GL_n(\HH)=\SL_n(\HH) \times \R^{>0}$ for $n= 1, 2$
and the isomorphism $K_2(\HH)\simeq H_2(\SL_1(\HH))$
from \cite[p. 287]{sah}.
\end{rem}

\section{Third homology of $\SL_2$ and the indecomposable $K_3$}

In this section we assume that $R$ is a commutative ring with many units,
unless it is mentioned. When a group $G$ acts on a module
$M$, we use the standard definition $M_G$ for  $H_0(G, M)$.
Consider the action of $\rr$ on $\SL_n$ defined by
\begin{gather*}
a. B:=\left(
\begin{array}{cc}
a & 0   \\
0 & 1
\end{array}
\right)
B
\left(
\begin{array}{cc}
a^{-1} & 0   \\
0       & 1
\end{array}
\right),
\end{gather*}
where $a \in \rr$ and $B \in \SL_n$.
This induces an action of $\rr$ on $H_i(\SL_n)$.
So by $H_i(\SL_n)_\rr$ we mean $H_0(\rr, H_i(\SL_n))$.

\begin{thm}\label{sl-in-sl}
Let $k$ be a field such that $1/2 \in k$.
\par {\rm (i)} $H_3(\SL_2,k)_\rr \arr H_3(\SL,k)$ is injective.
\par {\rm (ii)} If $R$ is an infinite field, then
$H_3(\SL_2,\zzz)_\rr \arr H_3(\SL,\zzz)$ is injective.
\par {\rm (iii)} If $R=\R$ or $R$ is an infinite field such that $\rr=\rr^2$, then
$H_3(\SL_2) \arr H_3(\SL)$ is injective.
\par {\rm (iv)} $H_3(\SL_2(\HH)) \arr H_3(\SL_3(\HH))$ is bijective.
\end{thm}
\begin{proof}
The part (iv) follows
from \ref{mir-elb4} and applying the K\"unneth
theorem to  $\GL_n(\HH)=\SL_n(\HH) \times \R^{>0}$, $n\ge 1$.

Since $H_3(\SL) \arr H_3(\GL)$ is
injective, to prove (i), (ii) and (iii), by \ref{mir-elb4} it
is sufficient to prove that
$H_3(\SL_2,k)_\rr \arr H_3(\GL_2,k)$,
$H_3(\SL_2,\zzz)_\rr \arr H_3(\GL_2,\zzz)$ and
$H_3(\SL_2) \arr H_3(\GL_2)$
are injective.

Set $A:=\zzz$ or $k$. From the map
$\gamma:\rr \times\SL_2 \arr \GL_2$, $(a, M) \mt aM$,
we obtain two short exact sequences
\begin{gather*}
1 \arr \mu_{2, R} \arr \rr \times\SL_2 \arr \im(\gamma) \arr 1,\\
1 \arr \im(\gamma) \arr \GL_2 \arr \rr/\rr^2 \arr 1.
\end{gather*}
Writing the Lyndon-Hochschild-Serre spectral sequence of the
above exact sequences and carrying out simple analysis, one gets
\begin{gather*}
H_3(\im(\gamma), A) \simeq H_3(\rr \times\SL_2, A),\ \
H_3(\im(\gamma),A)_{\rr/\rr^2}\simeq H_3(\GL_2, A).
\end{gather*}
Since the action of $\rr^2$ on $H_3(\im(\gamma), A)$
is trivial,
\begin{gather*}
H_3(\im(\gamma),A)_\rr
\simeq H_3(\GL_2, A).
\end{gather*}
These imply
\begin{gather*}
H_3(\GL_2, A) \simeq H_3(\rr \times\SL_2, A)_\rr.
\end{gather*}
Now the K\"unneth theorem implies that
$H_3(\SL_2,A)_\rr \arr H_3(\GL_2,A)$ is injective. This proves
parts (i) and (ii).

(iii) First let $\rr=\rr^2$. The map $\gamma$ induces the short exact sequence
\begin{gather*}
1 \arr \mu_{2, R} \arr \rr \times \SL_2 \arr \GL_2 \arr 1.
\end{gather*}
From the Lyndon-Hochschild-Serre spectral sequence of this
exact sequence one sees that $H_3(\inc): H_3(\SL_2) \arr H_3(\GL_2)$
has a kernel of order dividing $4$. To show that this kernel is trivial
we look at the spectral sequence induced by
$1 \arr \SL_2  \arr \GL_2 \arr \rr \arr 1$,
\[
{E'}_{p,q}^2=H_p(\rr, H_q(\SL_2)) \Rightarrow H_{p+q}(\GL_2).
\]
By \ref{sk1-k2} and the fact that the action
of $\rr$ on $H_i(\SL_2)$ is trivial , we get
the following ${E'}^2$- terms
\begin{gather*}
\begin{array}{cccccc}
\ast      & \ast               &                          &         &
           &  \\
H_3(\SL_2)& \ast               & \ast                     &         &
           &  \\
K_2^M(R)  &\rr \otimes K_2^M(R)& {E'}_{2,2}^2             & \ast    &
           &  \\
0         & 0                  & 0                        & 0       &
  0        &  \\
\z        & H_1(\rr)          & H_2(\rr)                & H_3(\rr)  &
H_4(\rr).&
\end{array}
\end{gather*}
Here ${E'}_{2,2}^2=H_2(\rr)\otimes K_2^M(R) \oplus
\tors(\mu_R,K_2^M(R))$ which is 2-divisible as $K_2^M(R)$ is uniquely
2-divisible. Hence
\[
H_3(\SL_2)/\im({d'}_{2,2}^2) \simeq {E'}_{0,3}^\infty \se H_3(\GL_2)
\]
which is induced
by $\SL_2 \hookrightarrow \GL_2$.
Thus $\im({d'}_{2,2}^2) \se \ker(H_3(\inc))$. This means that
$\im({d'}_{2,2}^2)$ is 2-divisible of order dividing $4$.
This is possible only if $\im({d'}_{2,2}^2)$ is trivial.

Now let $R=\R$. Consider the following exact sequences
\begin{gather*}
0 \arr \z/4\z \arr H_3(\SL_2(\R)) \arr H_3(\PSL_2(\R)) \arr 0,\\
0\arr  H_3(\PSL_2(\R)) \arr H_3(\PGL_2(\R)) \arr \z/2\z \arr 0
\end{gather*}
(see \cite[App. C, C.10, Thm. C.14]{parry-sah}).
In the first exact sequence $\z/4\z$ is mapped onto the subgroup of
order $4$ generated by $w:=\mtx {0} {1} {-1} {0}$
(see \cite[p. 207]{parry-sah}). Set
$\alpha:H_3(\SL_2(\R)) \arr H_3(\GL_2(\R)) $.
From the diagram
\[
\begin{array}{ccc}
H_3(\SL_2(\R))  &
\!\!\!
{\larr}
& \!\!\!  H_3(\GL_2(\R))
\\
\!\!\!\Big\downarrow\vcenter{%
\rlap{$\scriptstyle{}$}}
&       & \!\!\!\Big\downarrow\vcenter{%
\rlap{$\scriptstyle{}$}} \\
H_3(\PSL_2(\R))
& \!\!\!\larr & \!\!\! H_3(\PGL_2(\R))
\end{array}
\]
and the above exact sequences one
sees that $\ker(\alpha)$ is of order dividing $4$.
Here we describe the $E_2$-terms $E_{1,2}^2$ and $E_{2,2}^2$
of the spectral sequence
\begin{gather*}
E_{p,q}^2=H_p(\R^\ast, H_q(\SL_2(\R)) \Rightarrow H_{p+q}(\GL_2(\R)),
\end{gather*}
which is associated to
$1 \arr \SL_2(\R) \arr \GL_2(\R) \overset{\det}{\arr} \R^\ast \arr 1$.
It is well-known that 
\[
H_2(\SL_2(\R))\simeq K_2^M(\R)^\circ \oplus \z,
\]
where $K_2^M(\R)^\circ$ is the uniquely divisible part of
$K_2^M(\R)$. The action of $\R^\ast$ on $K_2^M(\R)^\circ$ is trivial
and its action on $\z$ is through multiplication by ${\rm sign}(r)$,
$r \in \R^\ast$ (see the proof of Prop. 2.15 in \cite[p. 288]{sah}).
Let $\bar{\z}$ be $\z$ with this new action of $\R^\ast$. Thus for $p=1,2$,
\[
E_{p,2}^2=H_p(\R^\ast) \otimes
K_2^M(\R)^\circ \oplus H_p(\R^\ast, \bar{\z}).
\]
It is not difficult to see that $H_1(\R^\ast, \bar{\z})=0$ and
$H_2(\R^\ast, \bar{\z})=\z/2\z$.
Now by an easy analysis of the above spectral sequence one sees that
$\ker(\alpha)$ is of order diving $2$. Since $w^2=-I_2 \in \GL_2(\R)$,
$\ker(\alpha)$, if not trivial, must be generated by
$x=[-I_2|-I_2|-I_2]$. But $\alpha(x)=[-I_2|-I_2|-I_2] \in H_3(\GL_2(\R))$
is non-trivial. Therefore $\ker(\alpha)=0$. Note that here
one has to use the fact that the action of $\R^\ast$ on
$H_3(\SL_2(\R))$ is trivial (see \cite[App. C.14]{parry-sah} and
\cite[2.10, p. 230]{dup-par-sah}). Therefore
$E_{0, 3}^2= H_3(\SL_2(\R))$.
\end{proof}

\begin{cor}\label{mir-sah}
Let $k$ be a field such that $1/2 \in k$.
\par
{\rm (i)} We have the split exact sequence
\begin{gather*}
0 \arr H_3(\SL_2,k)_\rr \arr H_3(\SL,k)
\arr K_3^M(R)\otimes k  \arr 0.
\end{gather*}
\par
{\rm (ii)} If $R$ is an infinite field, then we have the split exact sequence
\begin{gather*}
0 \arr H_3(\SL_2,\zzz)_\rr \arr H_3(\SL,\zzz)
\arr K_3^M(R)\otimes \zzz  \arr 0.
\end{gather*}
\par {\rm (iii)} If $R$ is an infinite field such that $\rr=\rr^2$, then
\begin{gather*}
0 \arr H_3(\SL_2) \arr H_3(\SL)
\arr K_3^M(R) \arr 0
\end{gather*}
is split exact.
\par {\rm (iv)}
We have the  split exact sequence
\[
0 \arr H_3(\SL_2(\R)) \arr H_3(\SL(\R))
\arr K_3^M(\R)^\circ \arr 0,
\]
where
$K_3^M(\R) \simeq \lan \{-1, -1, -1\} \ran \oplus K_3^M(\R)^\circ$.
\end{cor}
\begin{proof}
First we  prove (iv). The injectivity follows from \ref{sl-in-sl}.
From the diagram
\[
\begin{array}{ccccccccc}
\!\!\! 1  & \!\!\!\larr &\!\!\!  \SL_{2}(\R)  &
\!\!\!
{\larr}
& \!\!\!  \GL_{2}(\R) &  \!\!\! \larr
&\!\!\! \ \ \R^\ast
& \!\!\!\larr & \!\!\! 1 \\
&  &\!\!\!\Big\downarrow\vcenter{%
\rlap{$\scriptstyle{}$}}
&       & \!\!\!\Big\downarrow\vcenter{%
\rlap{$\scriptstyle{}$}}
&  &\!\!\!\Big\downarrow\vcenter{%
\rlap{$\scriptstyle{}$}}
&      &  \\
\!\!\! 1  & \!\!\!\larr & \!\!\! \SL(\R)
& \!\!\!\larr & \!\!\! \GL(\R)
& \!\!\!  \larr   & \!\!\! \ \   \R^\ast    &\!\!\!
\larr & \!\!\! 1.
\end{array}
\]
we obtain a map of spectral sequences
\begin{gather*}
\begin{array}{ccc}
E_{p, q}^2=H_p(\R^\ast, H_q(\SL_{2}(\R))) &
\Rightarrow & H_{p+q}(\GL_{2}(\R))\\
\Big\downarrow\vcenter{%
\rlap{$\scriptstyle{}$}} &      &
\Big\downarrow\vcenter{%
\rlap{$\scriptstyle{}$}}       \\
{E'}_{p, q}^2=H_p(\R^\ast, H_q(\SL(\R)))&
\Rightarrow & H_{p+q}(\GL(\R))
\end{array}
\end{gather*}
%
which give us a map of filtration
\begin{gather*}
\begin{array}{ccccccc}
0=F_{-1} & \se & F_0 & \se F_{1} \se & F_{2}  & \se & F_3=H_3(\GL_{2}(\R))\\
         &     & \downarrow & \downarrow   & \downarrow  &  & \downarrow \\
0=F_{-1}'& \se & F_0'& \se F_{1} \se & F_{2}' & \se & F_3'=H_3(\GL(\R)).
\end{array}
\end{gather*}
Since $H_3(\SL_{2}(\R)) \arr H_3(\GL_{2}(\R))$ is injective,
$F_0=E_{0,3}^\infty \simeq H_3(\SL_{2}(\R))$.
It is easy to see that $E_{p,1}^\infty = {E'}_{p,1 }^\infty=0$,
$F_0'={E'}_{0,3}^\infty \simeq H_3(\SL(\R))$ and
$E_{3,0}^\infty \simeq {E'}_{3,0 }^\infty$.
Since
\[
H_2(\SL_{2}(\R))=\z \oplus K_2^M(\R)^\circ \arr
\z/2\z \oplus K_2^M(\R)^\circ = H_2(\SL(\R))
\]
is surjective, $E_{2,2}^\infty \harr  {E'}_{2,2 }^\infty$ with
$\coker(E_{2,2}^\infty \arr  {E'}_{2,2 }^\infty) \simeq \z/2\z$ (see
the proof of \ref{sl-in-sl}(iii)). By an easy analysis of the above
filtration one gets the exact sequence
\[
0 \arr H_3(\SL(\R))/H_3(\SL_{2}(\R)) \arr
H_3(\GL(\R))/H_3(\GL_{2}(\R)) \arr \z/2\z \arr 0.
\]
Therefore $H_3(\SL(\R))/H_3(\SL_{2}(\R))\simeq K_3^M(\R)^\circ$.
A splitting map can be constructed using the composition
$K_3^M(\R)^\circ \arr H_3(\GL(\R)) \arr H_3(\SL(\R))$.

The proof of (i), (ii) and (iii) are similar. In the proof
of (iii) we need the homology stability $H_2(\SL_2)=H_2(\SL)$ and in
the proof of (i) and (ii) we need the isomorphism
\[
H_1(\rr, H_2(\SL_2,\zzz))\simeq H_1(\rr, H_2(\SL, \zzz)).
\]
To prove the latter, consider the exact sequence 
\[
1 \arr \rr^2 \arr \rr \arr \rr/\rr^2 \arr 1.
\]
This induces a map of Lyndon-Hochschild-Serre spectral sequences, with coefficients in
$H_2(\SL_2,\zzz)$ and $H_2(\SL,\zzz)$ respectively, 
which one easily  obtains the commutative
diagram
\begin{gather*}
\begin{array}{ccc}
 H_1(\rr^2, H_2(\SL_2,\zzz))_\rr & \overset{\simeq}{\larr} & H_1(\rr, H_2(\SL_2,\zzz))\\
\Big\downarrow\vcenter{%
\rlap{$\scriptstyle{}$}}
& & \Big\downarrow   \\
H_1(\rr^2, H_2(\SL,\zzz)) & \overset{\simeq}{\larr} & H_1(\rr, H_2(\SL,\zzz)).
\end{array}
\end{gather*}
The action of $\rr^2$ on $H_2(\SL_2,\zzz)$ is trivial, so 
\begin{gather*}
\begin{array}{ll}
H_1(\rr^2, H_2(\SL_2,\zzz))_\rr &
 \simeq (H_1(\rr^2,\zzz)\otimes H_2(\SL_2,\zzz))_\rr \\
 & \simeq  H_1(\rr^2,\zzz)\otimes H_2(\SL_2,\zzz)_\rr \\
 & \simeq H_1(\rr^2,\zzz)\otimes H_2(\SL,\zzz) \\
 & \simeq H_1(\rr^2, H_2(\SL, \zzz)).
\end{array}
\end{gather*}
Thus the left-hand column map in the above diagram is isomorphism. 
This implies  the
isomorphism of the right-hand column map.
\end{proof}
~
\begin{rem}
Let $R=\R$, $R=\HH$ or $R$ be an infinite field such that $\rr=\rr^2$.
Then $H_3(\SL_2)\arr H_3(\SL_3)$ is injective. This follows from Theorem
\ref{sl-in-sl}, and commutativity of the following diagram
\begin{gather*}
\begin{array}{ccc}
 H_3(\SL_2) &\!\!\!\!\!\!\!\!
\overset{}{\relbar\joinrel\!\relbar\joinrel\!\longrightarrow}  &
\!\!\!\!\!\!\!\!
H_3(\SL_3) \\
\ \ \ \ \ \  \searrow &
& \!\!\!\!\!\!\!\!\!\!\!\!\!\!\!\!\!\!\swarrow \\
& \!\!\!\!\!\!
H_3(\SL). &
\end{array}
\end{gather*}
This generalizes the main theorem of Sah in \cite[Thm. 3.0]{sah}.
\end{rem}
~

Let $K_3^M(R) \arr K_3(R)$ be the natural map from
the Milnor $K$-group to the Quillen $K$-group. Define
$K_3(R)^{\rm ind}:= \coker(K_3^M(R)\arr K_3(R))$.
This group is called the indecomposable part of $K_3(R)$.

\begin{prp}\label{mir-elb5}
Let $k$ be a field such that $1/2 \in k$.
\par {\rm (i)} $K_3(R)^{\rm ind}\otimes k \simeq H_3(\SL_2,k)_\rr.$
\par {\rm (ii)} If $R$ is an infinite field then
$K_3(R)^{\rm ind}\otimes \zzz \simeq H_3(\SL_2,\zzz)_\rr.$
\par {\rm (iii)} If $R=\R$, or if $R$ is an infinite field such that
$\rr=\rr^2$, then
$K_3(R)^{\rm ind} \simeq H_3(\SL_2)$.
\end{prp}
\begin{proof}
Let $A=\zzz, \z$ or $k$. By \ref{mir-sah} we have the commutative diagram
\[
\begin{array}{ccccccccc}
\!\!\! 0  & \!\!\!\larr &\!\!\! K_3^M(R)\otimes A&
\!\!\!
{\larr}
& \!\!\!  K_3(R) \otimes A &  \!\!\! \larr
&\!\!\! K_3(R)^{\ind} \otimes A
& \!\!\!\larr & \!\!\! 0 \\
&  &\!\!\!\Big\downarrow\vcenter{%
\rlap{$\scriptstyle{}$}}
&       & \!\!\!\Big\downarrow\vcenter{%
\rlap{$\scriptstyle{h_3}$}}
&  &\!\!\!\Big\downarrow{%
\rlap{$\scriptstyle{}$}}
&      &  \\
\!\!\! 0  & \!\!\!\larr & \!\!\!  K_3^M(R)\otimes A
& \!\!\!\larr & \!\!\! H_3(\SL,A)
& \!\!\!  \larr   & \!\!\!   H_3(\SL_2,A)_\rr &\!\!\!
\larr & \!\!\! 0.
\end{array}
\]
Here $h_3$ is the Hurewicz map $K_3(R)=\pi_3(B\SL^+) \arr H_3(\SL)$
and it
is surjective with two torsion kernel \cite[Prop. 2.5]{sah}.
In case $\rr=\rr^2$, $h_3$ is an isomorphism.
The snake lemma implies (i) and second part of (ii).
If $R=\R$, we look at the following commutative diagram
\[
\begin{array}{ccccccccc}
&  & K_3^M(\R) & {\larr} & K_3(\R)& \larr & K_3(\R)^{\ind} & \larr & 0 \\
&  &\Big\downarrow\vcenter{%
\rlap{$\scriptstyle{}$}}
&       & \Big\downarrow\vcenter{%
\rlap{$\scriptstyle{h_3}$}}
&  &\Big\downarrow{%
\rlap{$\scriptstyle{}$}}  &      &  \\
0 &\larr & K_3^M(\R)^\circ & \larr & H_3(\SL(\R))& \larr& H_3(\SL_2(\R))&
\larr &  0.
\end{array}
\]
The claim follows from the snake lemma using the fact that
$\ker(K_3(\R) \overset{h_3}{\arr} H_3(\SL(\R)))=\z/2\z$
\cite[2.17]{sah}.
\end{proof}

\begin{rem}
Theorem \ref{mir-elb5} generalizes theorem \cite[Thm. 4.1]{sah},
where three torsion is not treated.
\end{rem}

We can offer the following non-commutative version of the
above results.

\begin{prp}\label{h3-sl}
{\rm (i)} Let $R$ be a quaternion algebra. Then
\[
0 \arr H_3(\SL_2, \zzz)_\rr \arr
H_3(\SL, \zzz) \arr K_3^M(R) \otimes \zzz \arr 0
\]
is exact.
\par {\rm (ii)}
If $R$ is an Azumaya $R$-algebra, $R$ a commutative
local ring with an infinite residue field, then
\[
0 \arr H_3(\SL_2, \q)_\rr \arr H_3(\SL, \q) \arr
K_3^M(R) \otimes \q \arr 0
\]
is exact.
\end{prp}
\begin{proof}
(i) From the commutative diagram
\[
\begin{array}{ccccccccc}
\!\!\! 1  & \!\!\! \larr & \!\!\!  \SL_2  &
\!\!\! 
{\larr}
& \!\!\!  \GL_2 &  \!\!\! \larr
&\!\!\! \ \ K_1(R)
& \!\!\!\larr & \!\!\! 1 \\
&  &\!\!\!\Big\downarrow\vcenter{%
\rlap{$\scriptstyle{}$}}
&       & \!\!\!\Big\downarrow\vcenter{%
\rlap{$\scriptstyle{}$}}
&  &\!\!\!\Big\downarrow\vcenter{%
\rlap{$\scriptstyle{}$}}
&      &  \\
\!\!\! 1  & \!\!\!\larr & \!\!\! \SL
& \!\!\!\larr & \!\!\! \GL
& \!\!\!  \larr   & \!\!\! \ \   K_1(R)
& \!\!\!
\larr & \!\!\! 1
\end{array}
\]
we obtain a map of spectral sequences
\begin{gather*}
\begin{array}{ccc}
E_{p, q}^2=H_p(K_1(R), H_q(\SL_2, \zzz)) &
\Rightarrow & H_{p+q}(\GL_2, \zzz)           \\
\Big\downarrow\vcenter{%
\rlap{$\scriptstyle{}$}}
&   &
\Big\downarrow\vcenter{%
\rlap{$\scriptstyle{}$}}                \\
{E'}_{p, q}^2=H_p(K_1(R), H_q(\SL, \zzz))  &
\Rightarrow & H_{p+q}(\GL, \zzz)
\end{array}
\end{gather*}
Since the map $Z(\rr) \times \SL_2 \arr \GL_2$, $(a, B) \mt aB$,
has two torsion kernel and cokernel (use example \ref{k1-zr}),
$H_i(\SL_2, \zzz)_\rr \harr H_i(\GL_2, \zzz)$
(see the proof of \ref{sl-in-sl}(i)). 
By Lemma \ref{H-1}, $H_i(Z(\rr), \zzz) \harr H_i(\GL_2, \zzz)$ and it is easy to prove
the injectivity of $H_i(\SL, \zzz) \harr H_i(\GL, \zzz)$.
%
%
By an easy analysis of the above spectral sequences,
as in the proof of Cor. \ref{mir-sah},
we get the desired result. The proof of (ii) is similar.
\end{proof}

\begin{cor}
Let $D$ be a finite-dimensional $F$-division algebra.
Let 
\[
K_3^M(F,D):=\ker(K_3^M(F)\arr K_3^M(D)).
\]
Then we have the
following exact sequence
\[
0 \arr H_3(\SL_2(F), \q)_\fff \arr H_3(\SL_2(D), \q)_{D^\ast} \arr
K_3^M(F,D)\otimes \q \arr 0.
\]
\end{cor}
\begin{proof}
By Cor. 2.3 from \cite{hazrat2006},
$K_3(F) \otimes \q \simeq K_3(D)\otimes \q$.
Therefore \[
H_3(\SL(F), \q)\simeq H_3(\SL(D), \q)
\]
(see \cite[Thm. 2.5]{sah}). Now the claim
follows from Cor. \ref{mir-sah} and Prop. \ref{h3-sl}.
\end{proof}


\bigskip

\address{{\footnotesize
\
\\
Department of Pure Mathematics\\
Queen's University\\
Belfast  BT7 1NN \\
Northern Ireland\\
Email:\ b.mirzaii@qub.ac.uk
}}

\begin{thebibliography}{99}

\bibitem{arlettaz1990}
Arlettaz, D. A splitting result for the second
homology group of the general linear group.
Adams Memorial Symposium on Algebraic Topology, 1,
83--88, London Math. Soc. Lecture Note Ser., 175,
1992.

\bibitem{bass-tate}
Bass, H.; Tate, J. The Milnor ring of a global field.
Algebraic $K$-theory, II,
Lecture Notes in Math., Vol. 342, (1973), 349--446.

\bibitem{bor-yang}
Borel, A.; Yang, J. The rank conjecture for number fields.
Math. Res. Lett. {\bf 1} (1994), no. 6, 689--699.

\bibitem{bro}
Brown, K. S. Cohomology of groups.  Graduate
Texts in Mathematics, 87. Springer-Verlag, New York, 1994.

\bibitem{dennis1976}
Dennis, K.: In Search of New "Homology" Functors Having a
Close Relationship to K-theory, preprint (1976).

\bibitem{dup-par-sah}
Dupont, J. L., Parry, W., Sah, C.
Homology of classical Lie groups made discrete. II.
$H\sb 2,H\sb 3,$ and relations with scissors congruences.
J. Algebra {\bf 113} (1988), no. 1, 215--260.


\bibitem{elb}
Elbaz-Vincent, P. The indecomposable $K\sb 3$ of rings and homology of
$\SL_2$. J. Pure Appl. Algebra {\bf 132} (1998), no. 1, 27--71.

\bibitem{guin1989}
Guin, D. Homologie du groupe linéaire et $K$-théorie de
Milnor des anneaux. J. Algebra {\bf 123} (1989), no. 1, 27--59.

\bibitem{hazrat2006}
Hazrat, R. Reduced K-theory for Azumaya Algebras.
J. Algebra {\bf 305} (2006), 687--703.

\bibitem{lam2001}
Lam, T. Y. A first course in noncommutative rings. Second edition.
Graduate Texts in Mathematics, 131, Springer-Verlag, 2001.

\bibitem{mac}
Mac Lane, S. Homology, New York; Springer-Verlag,
Berlin-G\"ottingen-Heidelberg 1963.

\bibitem{mirzaii2003} Mirzaii, B. Homology stability for unitary
groups II. $K$-Theory {\bf 36} (2005), no. 3--4, 305--326.

\bibitem{mirzaii2005} Mirzaii, B. Homology of $\GL_n$:
injectivity conjecture for $\GL_4$. To appear in Math. Annalen.


\bibitem{nes-sus}
Nesterenko Yu. P., Suslin A. A. Homology of the general
linear group over a local ring, and Milnor's $K$-theory.
Math. USSR-Izv. {\bf 34} (1990), no. 1, 121--145.

\bibitem{parry-sah}
Parry, W., Sah, C. Third homology of $\SL(2,\R)$
made discrete.  J. Pure Appl. Algebra  {\bf 30}  (1983),
no. 2, 181--209.

\bibitem{sah1986}
Sah, C. Homology of classical Lie groups made discrete. I.
Stability theorems and Schur multipliers.
Comment. Math. Helv.  {\bf 61}  (1986),  no. 2, 308--347.

\bibitem{sah}
Sah, C. Homology of classical Lie groups made discrete.
III. J. Pure Appl. Algebra {\bf 56} (1989), no. 3, 269--312.

\bibitem{sus14}
Suslin, A. A. Homology of ${\rm GL}\sb{n}$, characteristic
classes and Milnor $K$-theory. Proc. Steklov Math. {\bf 3}
(1985), 207--225.


\bibitem{vdkallen1977}
Van der Kallen, W. The $K\sb{2}$ of rings with many units.
Ann. Sci. \'Ecole Norm. Sup. (4) {\bf 10}  (1977), no. 4, 473--515.


\bibitem{yag}
Yagunov, S. On the homology of $\GL_n$ and the higher
Pre-Bloch groups. Canad J. Math. {\bf 52} no. 6 (2000), 1310--1338.
\end{thebibliography}
\end{document}